\documentclass[12pt,a4paper,final]{iopart}

\usepackage{iopams}  
\usepackage{graphicx}
\usepackage[breaklinks=true,colorlinks=true,linkcolor=blue,urlcolor=blue,citecolor=blue]{hyperref}

\usepackage{amssymb}
\newcommand{\R}{{\mathbb R}}

\newcommand{\N}{{\mathbb N}}
\newcommand{\C}{{\mathbb C}}

\def\ter{\hfill \vrule width 5 pt height 7 pt depth - 2 pt\smallskip}

\newtheorem{thm}{Theorem}[section]
\newtheorem{cor}[thm]{Corollary}

\newtheorem{lem}[thm]{Lemma}
\newtheorem{pro}[thm]{Proposition}
\newtheorem{remark}[thm]{Remark}

\begin{document}

\title[Global continuation of bistable monotone waves]{Global continuation of monotone waves for  bistable delayed equations with unimodal nonlinearities}

\author{Sergei Trofimchuk}
\address{Instituto de Matem\'atica y Fisica,
Universidad de Talca, Casilla 747, Talca, Chile}
\ead{trofimch@inst-mat.utalca.cl}

\author[]{Vitaly Volpert}
\address{Institut Camille Jordan, UMR 5208 CNRS, University Lyon 1, 69622 Villeurbanne, France
}
\address{INRIA Team Dracula, INRIA Lyon La Doua, 69603 Villeurbanne, France
}
\address{RUDN University, ul. Miklukho-Maklaya 6, Moscow, 117198, Russia}
\ead{volpert@math.univ-lyon1.fr}

\begin{abstract} 
We study  the existence of monotone wavefronts for  a general family of bistable reaction-diffusion equations with delayed reaction term $g$. Differently from  
previous works, we do not assume the monotonicity of  $g(u,v)$  with respect to the delayed variable $v$
that  does not allow to apply  the  comparison techniques. Thus our proof is based on a variant of the Hale-Lin functional-analytic approach to heteroclinic solutions of functional differential equations where Lyapunov-Schmidt reduction is done in appropriate weighted spaces of $C^2$-smooth functions.  This method requires a detailed analysis of associated linear differential Fredholm operators and their formal adjoints.  For two different types of $v-$unimodal functions $g(u,v)$, we prove the existence of a maximal continuous family of bistable monotone wavefronts.. Depending on the type of unimodality 
(equivalently, on the sign of the wave speed),  two different scenarios can be observed for the bistable waves: 1) independently on the size of delay, each bistable wavefront is monotone; 2) wavefronts are monotone for moderate values of delays and can oscillate for large  delays.  
  \end{abstract}

\vspace{2pc}

\newpage

\section{Introduction and main results}
The main objects of investigation in this work  are  traveling front solutions for the delayed reaction-diffusion equation  
\begin{equation}\label{17ng}  \hspace{-15mm}
u_t(t,x) = u_{xx}(t,x)  + g(u(t,x), u(t-\tau,x)), \quad u \in \R,\ (t,x) \in \R^2, \ \tau \geq 0, 
\end{equation}
in the particular  case when the reaction term $g$ satisfies the following bistability condition: 

\vspace{1mm}

\noindent {\bf (B)}  Function $g$ is $C^{1,\gamma}$-continuous on some set  $(\alpha, \beta)^2 \subset \R^2$. On the interval $(\alpha, \beta)$,  equation  $g(u,u)=0$ has exactly three solutions $e_1<e_2<e_3$ such that $g_1(e_j,e_j)+ g_2(e_j,e_j)<0$ and  $g_1(e_j,e_j) <0$ for $j=1,3$ (in the paper, we use the notations $g_j(u_1,u_2)$, $j=1,2,$ for partial derivatives $\partial g(u_1,u_2)/\partial u_j$, $j=1,2$). 

\vspace{1mm}

We recall that classical solution $u(t,x) = \phi(x +ct)$  of (\ref{17ng}) is called a {\it bistable 
traveling front }(in the sequel, we  shorten this name to  the word `{\it wavefront}'  which will be used both for  the solution  $u(t,x) = \phi(x +ct)$  and for its profile  $\phi(s)$)  propagating with the velocity $c$, if $\phi$ is bounded $C^2$-smooth  function satisfying $\phi(-\infty) = e_1$ and
$\phi(+\infty) = e_3$. Wavefront is called monotone if $\phi'(t) \geq 0, \ t \in \R$.  Replacing in the above definition the boundary condition $\phi(+\infty) = e_3$ with  a weaker restriction  
$\liminf_{t \to +\infty} > e_2$, we define a classical solution called a semi-wavefront. 
It is clear that each wavefront $\phi$ to (\ref{17ng})  has to satisfy the following boundary value problem for delayed differential equation  
\begin{eqnarray} \label{twe2ng} \hspace{-15mm}
\phi''(t) - c\phi'(t) + g(\phi(t), \phi(t-c\tau)) =0,  \quad t \in \R, \quad \phi(-\infty)=e_1, \   \phi(+\infty)=e_3.  \end{eqnarray}
There are several particular forms of problem (\ref{twe2ng}) for which the existence of solutions is known. 
The simplest of them appears when $c\tau=0$:  problem  (\ref{twe2ng}) is then without delay and it is well understood  \cite{VVV}. In consequence,   we are interested only in non-stationary wavefronts  and will consider speed $c\not=0$.   Another well studied particular case of (\ref{twe2ng}) is when the nonlinearity $g(u,v)$ is non-decreasing in $v$ for each fixed $u$  \cite{BGHR,FW,GSW,LZ,MaWu,sch,SZ,WLR}. Indeed, this kind of monotonicity  allows  a successful application of  the maximum principle and comparison techniques.  

However, if the condition $g_2(u,v)\geq 0$ does not hold, not so much is known even about the existence of wavefronts to delayed reaction-diffusion equation (\ref{17ng}).\footnote[1]{And practically nothing is known about the uniqueness of bistable wavefronts in the non-monotone case, cf. \cite{ADG}. } In fact, we are aware about only two  such works,  \cite{ADG,TV1}, where special  cases of equation (\ref{17ng}) were analysed by means of the Leray-Schauder topological degree argument. In particular, the following  model 
\begin{equation}\label{17nM}
u_t(t,x) = u_{xx}(t,x)  - u(t,x) + f(u(t-\tau,x)), \quad u \geq 0,\ (t,x) \in \R^2, 
\end{equation}
with $C^{1,\gamma}$-continuous nonlinearity  $f:[0,+\infty)\to [0,+\infty)$ satisfying  $$f(0)=0=:e_1, \ f'(0) \in (0,1],\ f(e_2)-e_2=f(e_3)-e_3=0$$ 
has been  recently considered in \cite{ADG}.   
 Note that  bistable equation (\ref{17nM}) with the  unimodal birth function $f$ having only three fixed points, $e_1,e_2,e_3$, is broadly used in the mathematical ecology for  modelling systems exhibiting the Allee effect, cf. 
 \cite[Fig. 1d]{BY}.   Here, the unimodality of $f$ means that $f$ is hump-shaped, i.e.  it has a unique critical point, $\kappa$, and $0<e_2 < \kappa < e_3$. In \cite{GSW}, the above equation with unimodal $f$ was classified as Type D nonlinearity (see \cite[Fig. 4.4]{GSW}) and it was noted  in \cite[p. 5133]{GSW} that there has been no progress for Type D at the moment of the publication of \cite{GSW}.  In this regard, the recent contribution   \cite{ADG} by Alfaro, Ducrot and Giletti  presents a series of important existence results. 
Under some general bistability type assumptions on $f$ (which imply the positivity of the speed of propagation but are weaker than the unimodality restriction),  Alfaro {\it et al} 
 proved the existence of  semi-wavefront solutions to (\ref{17nM}) and established conditions sufficient for their either convergence or oscillation at $+\infty$.  In this paper, by giving a criterion for the existence of monotone wavefronts to (\ref{17ng}),  we provides an additional  insight into the interesting findings of \cite{ADG}.

\vspace{1mm}

Another type of bistable equation (\ref{17ng}) with unimodal nonlinearity was recently proposed in \cite{P1} in order to understand spatiotemporal dynamics of virus infection  spreading in tissues. The model equation of  \cite{P1} is of the form 
\begin{equation}\label{17nl}  \hspace{-15mm}
u_t(t,x) = u_{xx}(t,x)  + u(t,x)(1-u(t,x)-f(u(t-\tau,x))), \quad u \geq 0,\ (t,x) \in \R^2.
\end{equation}
It is assumed that  $f: \R_+\to (0,+\infty)$ is $C^{1,\gamma}$-continuous function  and that the equation $1-u=f(u)$
has exactly three positive solutions $0<e_1 <e_2<e_3$ on the interval $[0,1]$. In addition,    $f'(e_1) \geq 0$, 
$f'(e_3) > -1$ and  $f$ achieves the global maximum  at its unique critical point $\kappa \in (e_1,e_2)$.  See  Fig. 1 in \cite{P1}.  The recent work \cite{TV1} establishes that 
a simpler version of (\ref{17nl}) (with $e_1< 0 <e_2<e_3=1$)  has at least one monotone wavefront connecting the equilibria 0 and 1 for each fixed delay $\tau \geq 0$. 

\vspace{2mm}

The above mentioned biological models show the importance  of  studying the existence of wavefronts for equation (\ref{17ng}) with the reaction term  $g(u,v)$ which is not increasing in the second variable, but still has reasonably good (piece-wise monotone, with only two pieces of monotonicity) behavior  with respect to $v$ for each fixed $u$.  We will include  both models (\ref{17nM}) and (\ref{17nl}) in our general theory  by  considering two following  alternative unimodality assumptions:  

\vspace{2mm}

\noindent {\bf (U)} For each fixed $u\in (\alpha,\beta)$,  function $g(u,\cdot)$ has  
 a unique critical point $\kappa \in (e_1,e_2)$, independent on $u$ (hence, $g_2(u,\kappa)=0$) such that 
 $g_2(u,v)<0$ for $v\in (\alpha,\kappa)$ and   $g_2(u,v)>0$ for $v\in (\kappa, \beta)$.  
Furthermore, $g_1(u,v) <0$ for all $u \geq v$ such that $u \in [e_1, e_2)$, $v \in [e_1, \kappa]$ and 
$g(u,e_1)<0$ for all $u\in(e_2,\beta)$ while $g(u,e_1)>0$ for all $u\in(\alpha, e_1)$.  (The latter implies that $g(u,v) <0$ for all $u\geq v$, $u,v \in (e_1,e_2)$). 

\vspace{2mm}

\noindent {\bf (U$^*$)}  For each fixed $u\in (\alpha,\beta)$,  function $g(u,\cdot)$  has  
 a unique critical point $\kappa \in (e_2,e_3)$, independent on $u$  such that 
 $g_2(u,v)>0$ for $v\in (\alpha,\kappa)$ and   $g_2(u,v)<0$ for $v\in (\kappa, \beta)$.  
Furthermore, $g_1(u,v) <0$ for all $u \geq v;\ u,v \in [\kappa,e_3]$. In addition, $g(u,v)$ is `strongly' sub-tangential 
at $e_3$: $g_j(u,v) \geq g_j(e_3,e_3), \ j =1,2,  \ u \geq v, \ u, v \in [e_1,e_3]$. 

\vspace{2mm}

There is certain asymmetry in the strength of assumptions  {\bf (U)}  and {\bf (U$^*$)}: in particular, the sub-tangency requirement of {\bf (U$^*$)}   is used repeatedly in the proof of one of our main results, Theorem  \ref{main1b}.  Clearly, the `strong' sub-tangency condition is somewhat stronger than  the usual  sub-tangency requirement 
$g(u,v)\leq g_1(e_3,e_3)(u-e_3)+g_2(e_3,e_3)(v-e_3),$ $ u \geq v,$ $u,v \in [e_1,e_3]$. On the other hand, the form of sub-tangency given in  {\bf (U$^*$)} seems to be more friendly for applications. For instance, it is easy to see  that the reaction term in  (\ref{17nM}) satisfies {\bf (U$^*$)} if $f'(e_3) = \min\{f'(u), u \in [e_1,e_3]\}$ (note also that hypothesis {\bf (U)} holds for equation (\ref{17nl}) without additional restrictions on $f$).   Importantly, in Section \ref{TM} we show how a slightly weaker version of Theorem  \ref{main1b}, Theorem  \ref{main1c}, can be obtained without any kind of sub-tangency restriction  at $e_3$.

\vspace{1mm}

As we have mentioned, in this paper we consider only non-stationary wavefronts. In fact, 
it suffices to analyse the case of positive speeds, $c>0$, since  the  linear change of variables $\psi(-t)=e_1+e_3-\phi(t)$
transforms  problem (\ref{twe2ng}) under assumption {\bf (U)} and with the speed $c$ into  problem (\ref{twe2ng}) under assumption {\bf (U$^*$)} and with the speed $-c$, and vice versa (of course, modulo the sub-tangency condition and secondary monotonicity details).  In the next section, we are applying this trick in the case of  models (\ref{17nM}) and (\ref{17nl}).  Note also that if velocity $c$ is positive then traveling front is an expansion 
wave (since $\phi(x+ct)$ converges,  uniformly on compact sets, to the biggest steady state $e_3$ as $t\to +\infty$).
As we show, for the positivity of speed (for each fixed $\tau \geq 0$) it is enough to assume the inequality 

\vspace{2mm}

\noindent {\bf (I)} $\frak{I}:=  \int_{e_1}^{e_3}g(u,u)du >0$.   

\vspace{2mm}

Now, even if equation (\ref{17ng})  generally defines  a non-monotone evolutionary system,  we are interested in the existence of monotone wavefronts for it, cf. \cite{BNPR,FZ,FW,GTLMS,HW,TV1}. In the paper, such wavefronts will be obtained via deformation of the unique monotone wavefront of  equation (\ref{17ng})  considered with $\tau=0$.  The procedure 
of this continuous deformation requires from solutions of  (\ref{17ng})     the following monotonicity property (satisfied for  both considered biological models), cf.  \cite{GTLMS,TV1,VVV}:  

\vspace{2mm}

\noindent {\bf (M)}  Suppose that $u=\phi(x+ct),\ c >0$, is a non-decreasing wavefront 
connecting the steady states $e_1$ and $e_3$. Then $\phi'(t)>0, \ t \in \R$.   

\vspace{2mm}
 
Monotonicity of the  initial wavefront should be preserved during its continuous deformation.  It appears that it is easier to satisfy this requirement under assumption {\bf (U)} than under {\bf (U$^*$)}.  Indeed, as we will show in Lemmas \ref{cot}, \ref{2af+}, \ref{2af}, \ref{2aaf},  {\bf (U)} assures that each wavefront is strictly increasing at $\pm \infty$ and it is confined between the equilibria $e_1$ and $e_3$.  Contrary to this, if {\bf (U$^*$)} is assumed then it is easy to control monotonicity  at $-\infty$ but not at  $+\infty$ (asymptotic behaviour of monotone wavefronts at $+\infty$ is described in terms of zeros of the associated characteristic function $
\chi_-(z) = z^2-cz +a_- +b_- e^{-zc\tau},
$
where  coefficients $a_-, \ b_-$ are negative, see Section \ref{S4}). 
 A similar difficulty has occurred in \cite{GTLMS} during the continuous deformation of monostable monotone wavefronts.  In the cited work, it has been   shown that the monotone deformation of wavefronts can still be realised inside of some domain $\mathcal{D}$ of parameters $(\tau,c)$ described in continuation.   To define $\mathcal{D}$, we need the following result from  \cite[Lemma 1.1]{GTLMS} concerning the real zeros of $\chi_-(z)$:

\vspace{1mm}

\begin{pro} \label{PL} Given $a_-+ b_- <0, \ b_- < 0,$ there exists
$clin(\tau) \in (0, +\infty]$ such
that the characteristic equation $\chi_-(z)=0$, $c>0$, 
has  three real roots $\lambda_1\leq \lambda_2 <  0 < \lambda_3$
if and only if $c \leq clin(\tau)$.  If $clin(\tau)
$ is finite and $c=clin(\tau)$,  then $\chi_-(z)$ has a double zero $\lambda_1= \lambda_2<0$,
while for $c > clin(\tau)$ there does not exist any negative
root to  $\chi_-(z)=0$.  Moreover,  if $\lambda_j  \in \C$ is a
complex root of $\chi_-(z)=0$ for $c \in (0, clin(\tau)]$  then $\Re \lambda_j < \lambda_2$.

Furthermore,   $clin(\tau)= +\infty$ for all $\tau$ from some non-empty maximal interval $[0,\tau_\#]$
and $clin
(\tau)$ is strictly decreasing on $(\tau_\#, +\infty)$. In fact, \[  \hspace{-15mm} clin(\tau)=
\frac{\theta(a_-,b_-) +o(1)}{\tau}, \quad \tau \to
+\infty, \quad \mbox{ where} \
\theta(a_-,b_-):=
\sqrt{\frac{2\omega}{b_-}}e^{\omega/2},
\]
and $\omega$ is the unique negative root of
$
-2a_- = b_- e^{-\omega}(2+\omega).
$
\end{pro}
\begin{remark} \label{R1*} Suppose that $a_- <0$, then a straightforward analysis  shows that $\tau_\#>0$ can be determined as the unique real root of the equation $e|b_-|\tau e^{|a_-|\tau} =1$.  
For $\tau > \tau_\#$, 
the function $c=clin(\tau)$ can be defined implicitly by f the equation 
\begin{equation}\label{cL}  \hspace{-25mm}
\frak{A}(c,h):= \frac{2+\sqrt{c^2h^2+4+4|a_-|h^2}}{eh^2|b_-|} = \exp\left(\frac{2+2|a_-|h^2}{ch+\sqrt{c^2h^2+4+4|a_-|h^2}}\right) =:\frak{B}(c,h), 
\end{equation}
where $h = c\tau$.    Figures 2, 3 below present the graph of $c=clin(\tau)$ for $|a_-|=|b_-|=1$. 
\end{remark}
We define $\frak{ D}(a_-,b_-)$ as the set of non-negative parameters for which  $\chi_-(z), \ c >0,$ has exactly three real zeros (counting multiplicity). In the coordinates $(\tau,c)$, this domain takes the  next form
$$ \frak{ D}(a_-,b_-) = \{(\tau,c): \tau \geq 0, \  0 < c \leq clin(\tau)\} \subset \R_+^2. $$

We can  now state the first  main result of the paper: 
\begin{thm} \label{main1a}  Let assumptions {\bf (B)}, {\bf (I)},  {\bf (M)} and {\bf (U)} be satisfied. Then equation (\ref{17ng}) has a continuous family of strictly increasing  bistable wavefronts $u= \phi(x+c(\tau)t, \tau),$ $\tau \geq 0,$ propagating with the positive speed $c=c(\tau)$.  
\end{thm}
In Subsection \ref{toym}, we show (see Figure 3 below) that, under  assumptions of Theorem \ref{main1a}, it might happen that  $
c(\tau) > clin(\tau)$ for some positive values of $\tau$. 
Quite the contrary,  under assumption  {\bf (U$^*$)}, we need
the condition 
$$
(\tau,c(\tau)) \in Int\, \frak{ D}(g_1(e_3,e_3),g_2(e_3,e_3)), \ \mbox{where Int} \, \frak{ D} \ \mbox{denotes the interior of the domain}\  \frak{ D},  
$$ in order to realise monotone deformation of the initial wavefront: 
\begin{thm} \label{main1b}  Let assumptions {\bf (B)}, {\bf (I)},  {\bf (M)} and  {\bf (U$^*$)} be satisfied. Then there exists an extended real number $\tau_*> \tau_\#$ and a continuous function $c=c(\tau), \ \tau \in [0,\tau_*],$ such that 
equation (\ref{17ng}) has a continuous family of strictly increasing  bistable wavefronts $u= \phi(x+c(\tau)t, \tau),$ $\tau \leq \tau_*,$ propagating with the positive speed $c=c(\tau)$.  Moreover, 
$
(\tau,c(\tau)) \in \frak{ D}(g_1(e_3,e_3),g_2(e_3,e_3)),$ $c(\tau_*) = clin(\tau_*),
$
and $[0,\tau_*]$ is the maximal interval (containing $0$) for the existence of monotone wavefronts.  Furthermore, if $\tau_*$ is finite, then there is  a sequence of delays $\tau_j \to \tau_*$ such that 
 equation (\ref{17ng}) considered with $\tau=\tau_j$ has a wavefront propagating with speed $c_j$, $c_j \to  clin(\tau_*)$,  and oscillating around $e_3$. \end{thm}
 
In the next section, we apply Theorems \ref{main1a} and \ref{main1b}  to  models (\ref{17nM}) and (\ref{17nl}). In particular, we prove that condition {\bf (M)} is fulfilled for these equations.  Not only positive but also negative speeds of propagation are considered.  
In addition, in Subsection \ref{EX2}, we  state a somewhat weaker version of Theorem \ref{main1b}, Theorem \ref{main1c}. This result 
does not require any sub-tangency restriction from $g$.  In Subsection \ref{EX3}, 
 we are also illustrating our findings on an explicit example allowing a rather complete analytical and numerical analysis 
(this type of `toy models' was proposed in \cite{NRRP}, see also \cite{IGT,HK}).  In particular,   the computations done in Subsection \ref{EX3} suggest that  $c=c(\tau)$ is decreasing function of $\tau$ and that each monotone wavefront is unique (up to a translation). 

As in  \cite{GTLMS}, our proofs are based on the homotopy method and a variant of Hale-Lin
functional-analytic approach to the heteroclinic solutions \cite{HL}. In the bistable setting, this theory was developed  further by S.-N. Chow, X.-B. Lin, J. Mallet-Paret and W. Huang in  \cite{Chow,HW,HWU,MP}. In this theory, 
application of the Lyapunov-Schmidt reduction  requires a thorough analysis of the variational 
equations (and their adjoints) along the monotone wavefronts. Variational  equations are analysed in 
Section \ref{BU} (under assumption  {\bf (U)}) and Section \ref{S4} (under assumption  {\bf (U$^*$)}). 
The main conclusion of these sections concerns the existence of positive  (either on $\R$ or $\R_+$)  solutions $w_*(t)$ of the adjoint equations (Lemmas \ref{Lp} and \ref{Lp+}). Finally, Theorems \ref{main1a} and \ref{main1b} are proved in Section \ref{ThT}: to deal with the case when $w_*(t)$ can take negative values at some points $t <0$, we make
appropriate adjustments (expressed in terms of corrector functions) to the Lyapunov-Schmidt procedure. 
\section{Two biological  models and one illustrative example.}
\label{TM} 
In this section, we consider three different nonlinearities $g$ and, in each case, we apply the main results of the paper, Theorems \ref{main1a} and Theorems \ref{main1b}, to establish the existence of monotone (oscillating) wavefronts propagating with positive and negative speeds.  
\subsection{Mackey-Glass type model  (\ref{17nM}).}\label{EX1} \noindent  Assume that  $C^{1,\gamma}$-continuous unimodal function  $f: \R_+ \to \R_+$ satisfies 
\begin{itemize}
\item[a1)]  $f(0)=0=:e_1, \ f'(0) \in (0,1),\ f(e_2)-e_2=f(e_3)-e_3=0$. We also assume that equation $f(x)=x$ has only three solutions, $e_1,e_2,e_3$;
\item[a2)]  $f'(e_3) \leq f'(x),$ $x \in [0,e_3]$;
\item[a3)]  the unique critical point $\kappa$ of $f$ belongs to the interval $(e_2,e_3)$. 
\end{itemize}
Then $g(u,v) =-u +f(v)$ meets all restrictions of {\bf (B)}, {\bf (U$^*$)}. 
The wave profile equation for (\ref{17nM}) is 
\begin{eqnarray} \label{twe2an}  \hspace{-25mm}
\phi''(t) - c\phi'(t) - \phi(t) + f(\phi(t-c\tau)) =0,  \quad t \in \R, \quad \phi(-\infty)=e_1, \   \phi(+\infty)=e_3.  \end{eqnarray}
We claim that condition {\bf (M)} is satisfied in such a case. Indeed, let $\phi(t)$ be a profile of a bistable wave such that $\phi'(t) \geq 0$, $\phi(-\infty)=e_1$, $\phi(+\infty)=e_3$. Then Lemma \ref{2af*-} says that there exists a maximal interval $(-\infty, r)$ such that $\phi'(t) >0$ for all $t < r$. In addition, it holds that $\phi(r)\geq e_2$. Suppose that $r$ is finite, then $\phi'(r)=\phi''(r)=0$ so that $\phi(r) =  f(\phi(r-c\tau))$.  After differentiating (\ref{twe2an}), we also obtain that $\phi'''(r) = - f'(\phi(t-c\tau))\phi'(r-c\tau)$. Since $\phi'''(r) \geq 0$ and  $\phi'(r-c\tau)>0$, we find that $ f'(\phi(r-c\tau)) \leq 0$. Thus  $\phi(r-c\tau) \geq \kappa$ so that $\phi(r) \leq e_3 < f(\phi(r-c\tau))$, a contradiction. 
Hence, Theorem \ref{main1b} applies in such a case:
\begin{thm} \label{main1ba}  Let assumptions $a1), a2), a3)$ be satisfied together with  {\bf (I)} which here  reads as 
$$
\frak{P}:= \frac{1}{e_3-e_1}\int_{e_1}^{e_3}f(u)du -\frac{e_1+e_3}{2} >0. 
$$
Then all conclusions of Theorem \ref{main1b} are valid for equation  (\ref{17nM}).  \end{thm}

In Subsection \ref{EX3}, we present an explicit example showing that the Mackey-Glass type bistable models can have wavefronts oscillating around $e_3$. 

Next, in order to investigate the existence of monotone wavefronts for equation (\ref{twe2an}) when $\frak{P} <0$
we may apply, as it was suggested in the introduction, the change of variables $\psi(-t)=e_1+e_3-\phi(t)$. It transforms  the original equation  into  equation (\ref{twe2ng}) with new reaction term 
$\tilde g(u,v)=  e_1+e_3-u-f(e_1+e_3-v)$, steady states $e_1 < \tilde e_2 = e_1+e_3-e_2 < e_3$ and the critical 
point $\tilde \kappa = e_1+e_3-\kappa \in (e_1,\tilde e_2)$.  Moreover, it can be checked easily  that $\tilde g(u,v)$ satisfies  {\bf (B),  (U), (I)}  if  we assume conditions $a1), a3)$ and $\frak{P} <0$. 

\vspace{2mm}

Finally, let $\psi(t)$ be a wavefront  for the modified equation satisfying $\psi'(t) \geq 0$, $\psi(-\infty)=e_1$, $\psi(+\infty)=e_3$. Then Lemma \ref{cot} says that there exists a maximal interval $(-\infty, r)$ such that  $\psi'(t) >0$ for all $t < r$.  In addition, $\psi(r-|c|\tau)> \tilde \kappa$. Suppose that $r$ is finite, then $\psi'(r)=\psi''(r)=0$ so that $e_1+e_3-\psi(r) = f(e_1+e_3-\psi(r-|c|\tau))$, $\psi'''(r) = - f'(e_1+e_3-\psi(r-|c|\tau))\psi'(r-|c|\tau)$.  Since $
e_1+e_3-\psi(r-|c|\tau) < \kappa$, 
we conclude that $\psi'''(r)<0$, a contradiction. 
Hence, condition {\bf (M)} is satisfied by $\psi(t+|c|t)$ and an application of Theorem \ref{main1a}
leads to the following result. 
\begin{thm} \label{main1bb}  Assume conditions $a1), a3)$ as well as the inequality  $\frak{P} <0$. 
Then for each $\tau \geq 0$ equation (\ref{twe2an}) has a monotone wavefront  propagating with the negative speed $c(\tau)$ which depends continuously on the delay $\tau$. \end{thm}
\subsection{A model of virus infection spreading in tissues.}\label{EX2}  \noindent Following \cite{P1},  we consider reaction-diffusion equation  (\ref{17nl}) with the unimodal $C^{1,\gamma}$-continuous function $f: \R_+\to (0,+\infty)$ such that
\begin{itemize}
\item[b1)]  equation $1-u=f(u)$
has exactly three positive solutions $0<e_1 <e_2<e_3 <1$ on the interval $[0,1]$  and  $f'(e_1) \geq 0$,  
$f'(e_3) > -1$;
\item[b2)]  $f$ has a unique critical point $\kappa \in (e_1,e_2)$ where the global maximum of $f$ is achieved.
\end{itemize}
Then the function $g(u,v) =u(1-u -f(v))$ clearly satisfies the assumptions  {\bf (B)} (where $(\alpha, \beta)=(0,1)$) and  {\bf (U)}.  Next,  each bistable wavefront $\phi$ for model   (\ref{17nl}) solves the boundary problem 
\begin{eqnarray} \label{twe2an+}  \hspace{-25mm}
\phi''(t) - c\phi'(t) + \phi(t)(1-
\phi(t) - f(\phi(t-c\tau))) =0,  \  \phi(-\infty)=e_1 \geq 0, \   \phi(+\infty)=e_3. 
\end{eqnarray}
The assumption {\bf (M)} is also satisfied because of the following proposition.  
\begin{lem} \label{mono} Suppose that $\phi(t)$ satisfies (\ref{twe2an+}).  If  $\phi(t)$ is non-decreasing on some interval $(-\infty,s]$ and $\phi(t) \in [e_2,e_3)$ \ \ for 
$t \geq s$, 
then $\phi'(t)  >0, \ t \in \R$. \end{lem}
{\it Proof. }  Let $s' <s$ be a critical point for $\phi(t)$. Then $\phi'(s')=\phi''(s')=0$  so that 
$f(\phi(s'-c\tau)) = 1-\phi(s')$.  Since $\phi(s') \geq \phi(s'-c \tau)$ the latter equality implies that $\phi(s'-c\tau), \ \phi(s') >e_2$.  
Hence, $\phi(s-c\tau), \ \phi(s) >e_2$. Clearly, we may assume that $s = \sup\{r: \phi'(t) \geq 0, \ t \in (-\infty,r]\}$ and that there is $s_*\leq  s$ such that $\phi'(t) > 0$ for $t < s_*$, $\phi'(s_*)=0$. We have that either $s_*=s$ or $s_*<s$ and 
$\phi''(s_*)=0$.  In the latter case, $\phi(s_*-c\tau) >e_2$ and 
$$
\phi'''(s_*) = \phi(s_*)f'(\phi(s_*-c\tau))\phi'(s_*-c\tau)<0, 
$$
a contradiction.  Thus $\phi'(t)>0$ for all $t < s$ and, in addition, if $s$ is finite then $\phi''(s)<0$. 
Hence, if  $s$ is finite, then $\phi'(t) <0$ on some maximal interval $(s,S)$ (where $S$ is finite because of the condition $\phi(+\infty)=e_3$). Evidently, $\phi'(S)=0, \ \phi''(S)\geq 0$ so that $1-\phi(S) \leq f(\phi(S-c\tau))$. 
First, suppose that $S-c\tau \geq s$. Then $\phi(s) \geq \phi(S-c\tau)> \phi(S) \geq e_2$ implying that 
$$
1-\phi(S) \leq f(\phi(S-c\tau)) < f(\phi(S)), 
$$
a contradiction (since $f(x) \leq 1-x$ on $[e_2,e_3]$).  In consequence, $S-c\tau < s$ so that  
$$e_2 < \phi(s-c\tau) < \phi(S-c\tau) < \phi(s)> \phi(S)$$ yielding again a contradiction: 
$$
1-\phi(s) > f(\phi(s-c\tau)) > f(\phi(S-c\tau)) \geq 1-\phi(S).
$$
Thus $s=+\infty$ and $\phi'(t) >0$ for all $t \in \R$. 
\hfill \ter

An application of Theorem \ref{main1a}
allows us to extend the main existence result of \cite{TV1} on model (\ref{twe2an+}) considered under more realistic settings, cf. \cite{P1} :  
\begin{thm} \label{main1bbb}  Assume conditions $b1), b2)$ as well as condition {\bf (I)} which here  is equivalent to 
$$
\frak{J}:=  \int_{e_1}^{e_3}u(1-u-f(u))du >0. 
$$
Then for each $\tau \geq 0$ equation (\ref{twe2an+}) has a monotone wavefront  propagating with the positive speed $c(\tau)$ which depends continuously on the delay $\tau$. \end{thm}

Finally, supposing that $\frak{J}<0$, we will study the existence of wavefronts propagating with negative speeds 
(i.e. of the extinction waves).  Since we are going to invoke Theorem \ref{main1b}, this suggests the use of the transform $\psi(-t)=e_1+e_3-\phi(t)$.  The new reaction term has the form 
$
\hat g(u,v)= -(e_1+e_3- u) (1-e_1-e_3 + u- f(e_1+e_3-v))
$
and it is immediate to see that the 'strong' sub-tangency condition  of {\bf (U$^*$)} is not satisfied by 
$$
\hat g_1(u,v)=  1 -2 (e_1+e_3- u) - f(e_1+e_3-v), \quad g_2(u,v)= - (e_1+e_3- u)f'(e_1+e_3-v).  
$$

In this case, it is convenient to apply the following weaker version of Theorem \ref{main1b}:
\begin{thm} \label{main1c}  Let assumptions {\bf (B)}, {\bf (I)},  {\bf (M)} and  {\bf (U$^*$)} (except for the `strong' sub-tangency condition) be satisfied. Set 
\begin{equation}\label{abt} \hspace{-25mm}
\tilde a_- = \min \{g_1(u,v): u \geq v, \ u,v \in [e_1,e_3]\}, \ \tilde b_- = \min \{g_2(u,v): u \geq v, \ u,v \in [e_1,e_3]\}
\end{equation}
and let $\tilde \tau_\#>0$ be the unique real root of the equation $e|\tilde b_-|\tau e^{|\tilde a_-|\tau} =1$.  Then there exists an extended real number $\tau_*> \tilde \tau_\#$ and a continuous function $c=c(\tau), \ \tau \in [0,\tau_*]$, such that 
equation (\ref{17ng}) has a continuous family of strictly increasing  bistable wavefronts $u= \phi(x+c(\tau)t, \tau),$ $\tau \leq \tau_*,$ propagating with a positive speed $c=c(\tau)$.  Moreover, 
$
(\tau,c(\tau)) \in \frak{ D}(\tilde a_-,\tilde b_-)$ 
and the point $(\tau_*,c(\tau_*))$  belongs to the boundary of domain $\frak{ D}(\tilde a_-,\tilde b_-)$.   \end{thm}

\vspace{0.5mm}

\noindent \textsc{The strategy of the  proof of Theorem \ref{main1c}}.   The `strong' sub-tangency condition of  {\bf (U$^*$)} is invoked  four times in the proof of Theorem \ref{main1b}. In Remarks \ref{W1}, \ref{W2}, \ref{W3}, \ref{W4} below, we show how the exclusion of this condition changes the proof and the conclusion of  Theorem \ref{main1b}. 
\hfill \ter

\vspace{1mm}

Computing the parameters $\tilde a_-,  \tilde b_-$ and then applying  Theorem \ref{main1c} to the transformed version of equation (\ref{twe2an+}), we obtain   
 the following. 
\begin{thm} \label{main1bbc}  Assume conditions $b1), b2)$ and the  inequality  
$\frak{J} <0$. Set 
$$\tilde a_- = \min \{1-2u-f(u): u  \in [e_1,e_3]\}, \quad \tilde b_- = \min \{-uf'(u): u \in [e_1,e_3]\}$$ and define $\tilde \tau_\#>0$ as  in Theorem \ref{main1c}.  Then there is an extended real number $\tau_*> \tilde \tau_\#$ and a continuous function $c=c(\tau), \ \tau \in [0,\tau_*]$ such that 
model  (\ref{17nl}) has a continuous family of strictly increasing  bistable wavefronts $u= \phi(x+c(\tau)t, \tau),$ $\tau \leq \tau_*,$ propagating with negative speed $c=c(\tau)$.  Moreover, 
$
(\tau,|c(\tau)|) \in \frak{ D}(\tilde a_-,\tilde b_-)$ 
and the point $(\tau_*, |c(\tau_*)|)$  belongs to the boundary of the domain $\frak{ D}(\tilde a_-,\tilde b_-)$.     \end{thm}
{\it Proof. }   It is immediate to see that the function $\hat g(u,v)= -(e_1+e_3- u) (1-e_1-e_3 + u- f(e_1+e_3-v))$ meets conditions {\bf (B)}, {\bf (I)}  and  {\bf (U$^*$)} (except for the `strong' sub-tangency requirement). Now,  the validity of assumption {\bf (M)} for the transformed equation follows from a similar property of  the original model (\ref{twe2an+}) considered with $c<0$:  every its non-decreasing wavefront  $\phi(t)$
connecting  $e_1$ and $e_3$ satisfies the inequality $\phi'(t)>0, \ t \in \R$. In order to demonstrate this property, on the contrary, suppose that  the set $\mathcal S \subset \R$ of all critical points of $\phi$ is non-empty.   Then $\phi''(s)=0$ for each $s  \in \mathcal S$ and, consequently,
$1-\phi(s) = f(\phi(s-c\tau))$. Since $c <0$, this implies that $\phi(s) \in (e_1,e_2)$, $\phi(s-c\tau) \in (\kappa,e_2)$ so that 
$\sup {\mathcal S}=:s_0 \in \mathcal S$ is finite. After differentiating (\ref{twe2an+})  at $s_0$, we obtain the following contradiction  $0 \leq \phi'''(s_0) =\phi(s_0)f'(\phi(s_0-c\tau))\phi'(s_0-c\tau) <0$. Hence, $\mathcal S =\emptyset$ and  Theorem \ref{main1c}  can be applied to the equation with transformed reaction term $\hat g(u,v)$. 
\hfill \ter 

\subsection{A `toy' model.}\label{EX3} \label{toym} In this subsection, we are going to illustrate  results concerning the Mackey-Glass type model   by considering in (\ref{twe2an}) the 
following discontinuous nonlinearity  
\begin{eqnarray*}\label{KPPFt}
f(u)= \left\{\begin{array}{cc} 
pu,& u \in [0,\kappa), \\    
1 +q (u-1), & u \geq   \kappa,\end{array}\right. 
\end{eqnarray*}
with $\kappa, p \in (0,1), q <0$ and $e_1=0, \ e_3 =1$, see Figure 1.

 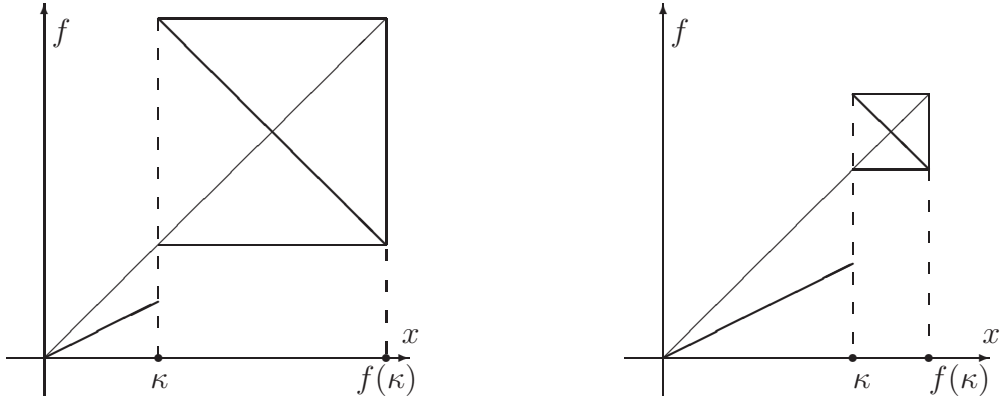
\begin{figure}[h]\label{F2}
 \setlength{\unitlength}{1cm}
\begin{picture}(8,5)(-1.5,-0.3)
\put(-0.5,0){\vector(1,0){5.3}}
\put(4.1,-0.4){$ f(\kappa)$}
\put(4.7,0.2){$x$}
\put(0.1,4.2){$f$}
\put(1.4,-0.4){$\kappa$}
\put(1.5,0){\circle*{0.1}}
\put(4.5,0){\circle*{0.1}}
\put(0,0){\line(1,1){4.5}}
\multiput(4.5,1.85)(0,-0.41){5}
{\line(0,-1){0.15}}
\multiput(1.5,4.5)(0,-0.395){12}
{\line(0,-1){0.15}}
\put(0,-0.5){\vector(0,1){5.2}}
\put(1.5,4.5){\line(1,0){3}}
\put(4.5,4.5){\line(0,-1){3}}
\put(4.5,1.5){\line(-1,0){3}}
\thicklines
\put(0,0){\line(2,1){1.5}}
\put(1.5,4.5){\line(1,-1){3}}
\end{picture}
\begin{picture}(1,5)(-1.5,-0.3)
\put(-0.5,0){\vector(1,0){4.8}}
\put(3.5,-0.4){$ f(\kappa)$}
\put(4.2,0.2){$x$}
\put(0.1,4.2){$f$}
\put(2.5,-0.4){$\kappa$}
\put(2.5,0){\circle*{0.1}}
\put(3.5,0){\circle*{0.1}}
\put(0,0){\line(1,1){3.5}}
\multiput(3.5,2.5)(0,-0.41){6}
{\line(0,-1){0.15}}
\multiput(2.5,3.5)(0,-0.395){9}
{\line(0,-1){0.15}}
\put(0,-0.5){\vector(0,1){5.2}}
\put(2.5,3.5){\line(1,0){1}}
\put(3.5,3.5){\line(0,-1){1}}
\put(3.5,2.5){\line(-1,0){1}}
\thicklines
\put(0,0){\line(2,1){2.5}}
\put(2.5,3.5){\line(1,-1){1}}
\end{picture}
\caption{\hspace{0cm} Graph of the unimodal birth function $f$: the cases $c > 0$ (on the left) 
and $c < 0$ (on the right).}
\end{figure}
First, we assume the inequality ${\bf (I)}$:  $\int_{0}^{1}(-u +f(u))du >0$, see Figure 1,  left.   This condition amounts to 
\begin{equation}\label{kla}
k_*:=\kappa \left(1+ \sqrt{\frac{1-p}{1-q}}\right) <1. 
\end{equation}
Let $\phi$ be a profile of a bistable wave normalised by the condition $\phi(-c\tau)=\kappa$. Clearly, 
$\phi$ is a positive solution of the linear equation
\begin{equation}\label{MGDtp}
\phi''(t)-c\phi'(t)-\phi + p\phi(t-c\tau)=0, \quad t <0. 
\end{equation}
The characteristic equation for (\ref{MGDtp}) is
\begin{equation}\label{MGDtc}
\lambda^2 -c\lambda - 1 + pe^{-c\tau \lambda}= 0,
\end{equation}
and it has a unique positive real root $\mu_1=\mu_1(c,\tau)$,  see also Lemma \ref{lc2} below. Thus
$$
\phi(t)=\kappa e^{\mu_1(t+c\tau)}, \quad t\leq 0.
$$
Hence, if $t>0$ and  $\phi(t) \geq \kappa$ for all $t\geq 0$ (this requirement is automatically satisfied for each monotone bistable wave), then $\phi(t)$ for $t > 0$ satisfies   the equation
$$
\phi''(t)-c\phi'(t)-\phi +1 +q (\phi(t-c\tau)-1)=0.
$$
The change of variables $\psi=\phi-1$ transforms this equation into
\begin{equation}\label{MGDtp2}
\psi''(t)-c\psi'(t)-\psi + q\psi(t-c\tau)=0.
\end{equation}
We also have that 
\begin{equation}\label{IvP}  \hspace{-15mm}
\psi(t) = \kappa e^{\mu_1(t+c\tau)}-1,\  t \in [-c\tau,0], \quad 
\psi(0)= \kappa e^{\mu_1c\tau}-1, \quad \psi'(0) = \kappa \mu_1 e^{\mu_1c\tau}.
\end{equation}
Applying the Laplace transform $(L\psi)(z)= \int\limits_0^{\infty}e^{-zt}\psi(t)dt$ to equation (\ref{MGDtp2}), we get
$$
\chi(z)(L\psi)(z)=\psi'(0)+z\psi(0)-c\psi(0)-qe^{-zc\tau}\int\limits_{-c\tau}^{0}e^{-zt}\psi(t)dt.
$$
Here $\chi(z)= z^2 -cz - 1 +qe^{-c\tau z}$ has a unique positive zero $\lambda_1$, see Lemma \ref{lc3}. Furthermore, we will assume that  the parameters $c, \tau, q$ are such that  $\lambda_1$ is the only zero of $\chi(z)$ on the closed right half-plane. Then the stable manifold of the zero equilibrium to (\ref{MGDtp2}) has codimension 1 and the solution of 
 initial value problem (\ref{IvP}) for this equation belongs to the stable manifold if and only if the projection of the initial function on the unstable manifold is zero, i.e. if  and only if 
$$
\kappa \mu_1 e^{\mu_1c\tau} + (\lambda_1-c)(\kappa e^{\mu_1c\tau}-1)-qe^{-\lambda_1 c\tau}\int\limits_{-c\tau}^{0}e^{-\lambda_1 t}\psi(t)dt = 0.
$$
After an integration, this gives 
$$
\kappa \left(\mu_1 e^{\mu_1c\tau} + (\lambda_1-c) e^{\mu_1c\tau}-q \frac{e^{(\mu_1-\lambda_1)c\tau}-1}{\mu_1-\lambda_1}\right) = (\lambda_1-c)+\frac{q}{\lambda_1}
(e^{-\lambda_1 c\tau}-1). 
$$
Since $\mu_1$ and $\lambda_1$ are solutions of equations (\ref{MGDtc}) and $\chi(z)= 0$,  respectively, the last equation simplifies to 
\begin{equation}\label{theta}
\kappa=\frac{1-q}{p-q}\left(1-\frac{\mu_1(c)}{\lambda_1(c)}\right) =:K(c) >0, \quad c \geq 0. 
\end{equation}
Being simple zeros, $\lambda_1(c)$ and $\mu_1(c)$ are positive continuous functions of $c$ and clearly  $K(+\infty)= 0$. In addition, due to (\ref{kla}), 
$$
K(0)= \frac{1-q}{p-q} \left(1-\frac{\mu_1(0)}{\lambda_1(0)}\right) = \frac{1-q}{p-q} \left(1-\sqrt{\frac{1-p}{1-q}}\right) =  \left(1+ \sqrt{\frac{1-p}{1-q}}\right)^{-1} > \kappa. 
$$
Thus for each delay $\tau \geq 0$ there exists at least one speed $c >0$ such that  equation (\ref{theta}) is satisfied. 
In fact, the next result shows that such $c$ is actually unique (and therefore,  for each fixed $\tau$,   bistable wavefront of the `toy' version of  (\ref{twe2an}) is unique up to  translation). 
\begin{lem} 
It holds that $(\mu_1(c)/\lambda_1(c))' >0$ for all $c >0$. 
\end{lem}
{\it Proof. }  For $c >0$, set $\epsilon = c^{-2}$ and observe that functions $\mu(\epsilon) := c\mu_1(c)$ and 
$\lambda(\epsilon) := c\lambda_1(c)$ satisfy the equations
$
\epsilon z^2 -z -1 + pe^{-z\tau} =0, \ \epsilon z^2 -z -1 + qe^{-z\tau} =0
$, respectively.  Clearly, the lemma statement amounts to $(\mu(\epsilon)/\lambda(\epsilon))' \not=0$.  
So, on the contrary, suppose  that the latter derivative is equal to $0$ at some point $\epsilon_0$. Then 
$\mu'(\epsilon_0)\lambda(\epsilon_0) = \mu(\epsilon_0)\lambda'(\epsilon_0)$.  Set $\lambda_0=  \lambda(\epsilon_0)$,  $\mu_0=  \mu(\epsilon_0)$. Since 
$$
\lambda'(\epsilon_0)= - \frac{\lambda^2_0}{2\epsilon_0\lambda_0-1+\tau(\epsilon_0\lambda_0^2-\lambda_0-1)}, \ \mu'(\epsilon_0)= - \frac{\mu^2_0}{2\epsilon_0\mu_0-1+\tau(\epsilon_0\mu_0^2-\mu_0-1)}, 
$$
the equality $\mu'(\epsilon_0)\lambda(\epsilon_0) = \mu(\epsilon_0)\lambda'(\epsilon_0)$ is equivalent to 
$$
 \frac{2\epsilon_0\lambda_0-1+\tau(\epsilon_0\lambda_0^2-\lambda_0-1)}{\lambda_0}=  \frac{2\epsilon_0\mu_0-1+\tau(\epsilon_0\mu_0^2-\mu_0-1)}{\mu_0}, 
$$
which can be simplified to the following contradictory relations
$$
- (1+\tau) \left(\frac 1\lambda_0 - \frac 1\mu _0\right) = \tau\epsilon_0 (\mu_0 - \lambda_0), \quad 0> -(1+\tau)= \tau\epsilon_0\mu_0\lambda_0 >0. \hspace{3cm} \ter
$$
\begin{figure}[h]
\centering \fbox{\includegraphics[width=10cm]{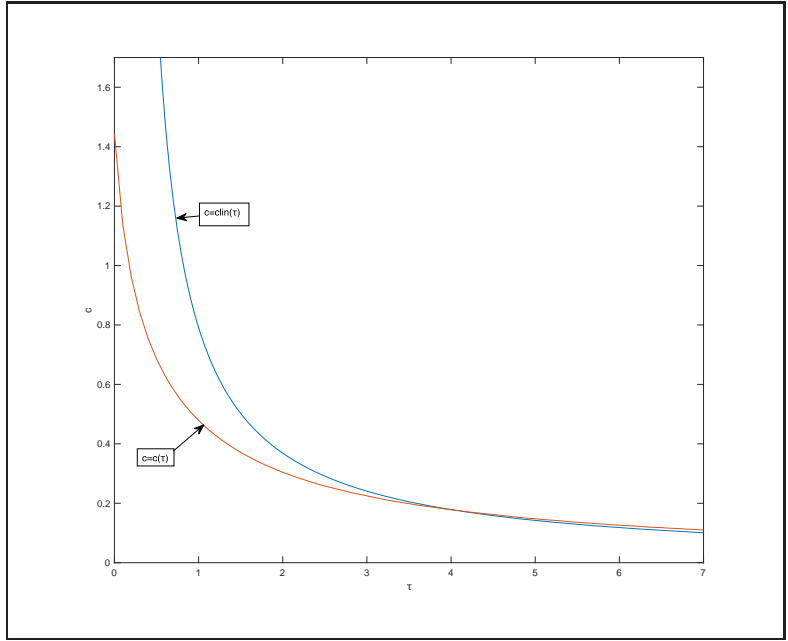}}
\caption{\hspace{-0.5cm} Domain $\mathcal{D}(-1,-1)$ and the curve $c=c(\tau)$ for the
`toy' model with $\kappa=1/3, \ p =1/2,\ q =-1$.}
\end{figure}
\noindent For numerics,  we take $\kappa=1/3, \ p =1/2, q =-1$ as shown on Figure 1 (left). 
Then for small values of delays ($0 \leq \tau \leq 4.04\dots$), our `toy' equation  has monotone bistable waves 
propagating with speeds $c=c(\tau)$, see Figure 2. However, for bigger 
delays (i.e. for $\tau > 4.04\dots$) the bistable wave profile $\phi(t)$ oscillates around the equilibrium $1$ at $+\infty$.  Figure 2 also suggests  that $c(\tau)$ is a decreasing function of $\tau$. 

If we now suppose that $k_*>1$ (i.e.  $ \int_{0}^{1}(-u +f(u))du <0$, see Figure 1, right),  then the propagation speed must  be negative, $c <0$. Using the same notations $\lambda_1(|c|), \mu_1(|c|)$ and applying the Laplace transform approach again, we find that for every delay $\tau\geq 0$ there exists a unique monotone bistable wave propagating with the speed $c=c(\tau)<0$ which can be determined from the equation
$$
1- \kappa=\frac{1-p}{p-q}\left(\frac{\lambda_1(|c|)}{\mu_1(|c|)}-1\right)  >0, \quad c < 0. 
$$
For $t \geq 0$, the explicit form of the unique profile normalised by the condition $\phi(-c\tau)=\kappa$ is given by 
$
\phi(t) = 1 - (1-\kappa)e^{-\lambda_1(|c|)(t+c\tau)}. 
$ 
\begin{figure}[h]
\centering \fbox{\includegraphics[width=10cm]{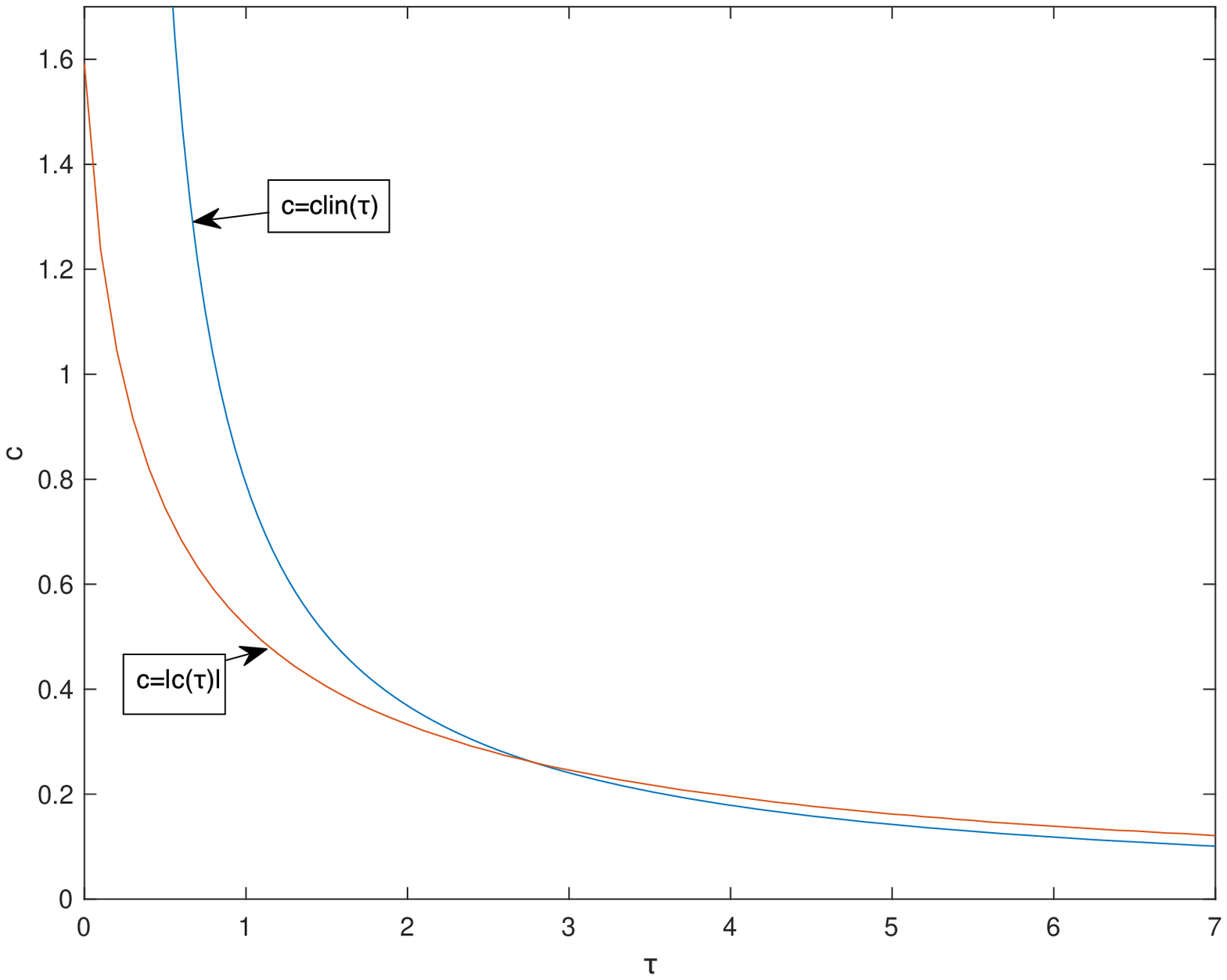}}
\caption{\hspace{-0.5cm} Domain $\mathcal{D}(-1,-1)$ and the curve $c=|c(\tau)|$ for the
`toy' model with $\kappa=0.9, \ p =1/2,\ q =-1$.}
\end{figure}
For numerical calculations in this case, we take $\kappa=0.9, \ p =1/2,\ q =-1$ as shown on Figure 1, right. 
Figure 3 also suggests that $|c(\tau)|$ is a decreasing function of $\tau$ and shows that the inequality $|c(\tau)| > clin(\tau)$  may happen for certain values of $\tau$ (unlike the case when $c>0$, the latter does not affect the monotonicity property of the profile $\phi$). 
\section{Variational equation along the monotone bistable wave under assumption {\bf (U)}.}\label{BU}
Let $\phi(t)$ be a solution of problem (\ref{twe2ng})  with $c>0$ satisfying $\phi(t) \geq e_1$ for all $t \in \R$.  
Without restricting generality, we can assume that  $\phi(-c\tau)= \kappa$ and that $\phi(t) < \kappa$ for $t < -c\tau$. 
The variational equation along $\phi(t)$ is of the form ${\mathcal D}\psi(t)=0$, where
$$
{\mathcal D}\psi(t):=\psi''(t)- c \psi'(t) + a(t)\psi(t) + b(t)\psi(t-c\tau), 
$$
and
$$
a(t):= g_1(\phi(t), \phi(t-c\tau)), \quad b(t) =  g_2(\phi(t), \phi(t-c\tau)). 
$$
Clearly, ${\mathcal D}\phi'(t)=0$. In view of assumptions   {\bf  (B),  (U)},  we have that
$
b(0)=0;$ \ $ a(t) <0$ for $t < -c\tau$; \ $b(t) <0$ for  $t < 0$  and $b(t) > 0$ for  $t > 0$;
$$
a_-:=a(-\infty)= g_1(e_1, e_1) < 0, \quad b_-:=b(-\infty) = g_2(e_1, e_1) < 0; 
$$
$$
 a_+:=a(+\infty)= g_1(e_3, e_3) { <0}, \quad b_+:= b(+\infty) = g_2(e_3, e_3) >0, \quad a_++ b_+ <0.
$$
Hence, the variational equation is asymptotically autonomous and the limiting autonomous equations at $\pm\infty$
have the  characteristic functions
$
\chi_\pm(z) = z^2-cz +a_\pm +b_\pm e^{-zc\tau}. 
$
It is easy to see that $\chi_+(z)$ always has exactly two real roots (we will denote them by $\mu_2 < 0 < \mu_1$) and that 
$\chi_-(z)$ always has exactly one positive root (we will use the notation $\lambda_1$ for it).   Some further information 
about zeros of $\chi_\pm(z)$ can be found in the next two lemmas. 
\begin{lem} \label{lc2} The zeros $\mu_1, \mu_2$ are simple. Moreover, they are unique zeros of $\chi_+(z)$ in the half-plane  $\{\Re z \geq \mu_2\}$.  
\end{lem}
{\it Proof. }  Since $\chi_+(x) <0$ for all $x \in (\mu_1,\mu_2)$ and $\chi_+''(x) >0, \ x \in \R$, the equalities 
$\chi_+'(\mu_j) =0,$  $j =1,2,$ are excluded. Thus $\mu_1, \mu_2$ are simple zeros.  Next, let $z_j^{\pm}, j =1,2,$ denote 
the real zeros of the polynomial  $z^2-cz +a_\pm$.Then $z_1^+< \mu_2<0< \mu_1 < z_2^+$.   If $w$ is a complex zero of $\chi_+(z)$, it holds that 
$$
|\Re w-z_1^+||\Re w-z_2^+|<|w-z_1^+||w-z_2^+| = b_+e^{-\Re w c\tau}, 
$$
so that, for each complex zero $w$ with $\Re w \in [\mu_1,\mu_2]$
$$
0< b_+e^{-\Re w c\tau}  - |\Re w-z_1^+||\Re w-z_2^+| = b_+e^{-\Re w c\tau}  - (\Re w-z_1^+)(z_2^+- \Re w) = \chi_+(\Re w),
$$
contradicting to  the inequality $\chi_+(x) <0, \ x \in (\mu_1,\mu_2)$.

Similarly, for each $z=iy, \ y \in \R, $ it holds
$$
|b_+e^{-z c\tau}| = b_+ <  -a_+= |\Re iy -z_1^+||\Re iy -z_2^+|\leq |z-z_1^+||z-z_2^+| = |z^2-cz+a_+|.  
$$
As a consequence, by a standard argument invoking the Rouch\'e theorem, the numbers of roots of $z^2-cz+a_+$ and $z^2-cz+a_+ + b_+e^{-z c\tau}$ on the half-plane $\{\Re z \geq 0\}$ coincide (due to decaying nature of $b_+e^{-x c\tau}$ for $x >0$).  
\hfill \ter
\begin{lem} \label{lc3} $\lambda_1$ is simple and dominating zero  of $\chi_-(z)$: every other root $\lambda_j$ of the equation $\chi_-(z)=0$ satisfies $\Re \lambda_j < \lambda_1$.  If $\lambda_2,\ \Im \lambda_2 \geq 0,$ is a root with the biggest 
real part $\Re \lambda_2 < \lambda_1$, then  $\lambda_2$ is a unique root of $\chi_-(z)=0$ with these properties belonging to the upper half-plane. Moreover, $\lambda_2$ is either  real negative root of the maximal multiplicity 2, or it is a simple complex root. 
\end{lem}
{\it Proof. }  Clearly, $\chi_-'(\lambda_1)>0$ and therefore the multiplicity of $\lambda_1$ is equal to 1. Next,  for each $z$ with $x=\Re z>\lambda_1$, we have that 
$$
|z^2-cz+a_-| = |z-z_1^-||z-z_2^-| \geq (x-z_1^-)(x-z_2^-)= x^2-cx+a_- > |b_-|e^{-c\tau x} = |b_-e^{-c\tau z}|, 
$$
so that every zero $\lambda_j$ of the characteristic function should satisfy $\Re \lambda_j \leq \lambda_1$.  
Now, if  $\Re z = \lambda_1,$  $\Im z \not=0$, then $|z-z_1^-||z-z_2^-| > (x-z_1^-)(x-z_2^-)$ so that again 
$|z^2-cz+a_-|  >  |e^{-c\tau z}|.$

Next, if $\lambda_j = \alpha +i \beta_j$, $0\leq \beta_2< \beta_3, \ j=2,3$, are 
two zero of $\chi_-(z), $ then 
$$
|\lambda_2^2-c\lambda_2+a_-|  =  |b_-|e^{-c\tau  \alpha} = |\lambda_3^2-c\lambda_3+a_-|. 
$$
On the other hand, 
$$
|\lambda_2^2-c\lambda_2+a_-|  = |\lambda_2-z_1^-||\lambda_2-z_2^-|  < |\lambda_3-z_1^-||\lambda_3-z_2^-| = 
|\lambda_3^2-c\lambda_3+a_-|,
$$
a contradiction which proves the uniqueness of $\lambda_2$.   Let us suppose now that the complex zero 
$\lambda_2$ is multiple. Then $\chi_-(\lambda_2)= \chi'_-(\lambda_2)=0$ that implies 
$\lambda_2^2-(c-2/(c\tau))\lambda_2 + a_- -1/\tau=0$.  Since the latter quadratic equation has only real roots, we 
get a contradiction.  
Finally,  inequality $a_- <0$ implies that the system of equations $\chi_-(\lambda_2)= \chi'_-(\lambda_2)=\chi''_-(\lambda_2)=0$ is incompatible.   In this way, the multiplicity of $\lambda_2$ cannot exceed two. 
\hfill \ter

Equation  ${\mathcal D}\psi(t)=0$ can be written as the
system
\begin{equation}\label{asy}
 v'(t)  = w(t), \
  w'(t) = -a(t)v(t)+ c w(t) - b(t)v(t-c\tau),
\end{equation}
or shortly as ${\frak F}_c(v,w)=0$, where \[ {\frak F}_c (v,w)(t)
= (v'(t)-w(t), w'(t)+  a(t)v(t)- c w(t) + b(t)v(t-c\tau)).
\]
System (\ref{asy}) possesses exponential dichotomy  at $+\infty$ and  shifted exponential dichotomy with exponents $\alpha_1:=\lambda_1-1.5\delta < \lambda_1-0.5\delta=: \beta_1$ (for $\delta >0$ small enough to satisfy 
$\Re \lambda_2 < \lambda_1-1.5\delta$) at $-\infty$, see \cite{HL} for the definition of these dichotomies. 
As a consequence, each solution of (\ref{asy})  converging to $0$ at $+\infty$, has an exponential rate of decay. Since for each bistable wavefront $\phi$ we have that $\phi''(\pm\infty)=\phi'(\pm\infty)=0$, we conclude that  $\phi', \phi''$ converge to $0$ at $+\infty$ with the exponential rate. More precise asymptotic formulas are given 
in Lemma \ref{2af+}.  To deal with the problem of super-exponentially small solutions in the proof of Lemma \ref{2af+}, we first establish the positivity of $e_3-\phi(t)$ for all $t\in \R$ : 
\begin{lem} \label{cot} Assume  {\bf (U)}  and suppose that $\phi(t) \geq e_1, \ t \in \R$, $\phi(-c\tau)= \kappa$,  is a bistable  wavefront for equation (\ref{twe2ng}) propagating with the speed $c>0$.  Then $\phi(t)  < e_3, \ t \in \R$.  If $\psi(t)$, $\psi(-\infty)=e_1,$   is a non-constant solution of equation (\ref{twe2ng}) which  is non-decreasing on the maximal interval $(-\infty, s]$ 
then $\psi(s-c\tau) > \kappa$. Consequently, if  $\psi(t)$ is normalized by the relation $\psi(-c\tau) =\kappa$ then  $s >0$.  Moreover, the normalized solution $\psi(t)$ satisfies $\psi''(t) >0$ and $\psi'(t)> 0$ for all  $t\leq 0$. 
\end{lem}
{\it Proof. }  (a) Indeed, otherwise $\phi(t)$ reaches its absolute maximum on $\R$ at some point $s_2$, where 
$\phi'(s_2)=0, \ \phi''(s_2)\leq 0, \  \phi(s_2) \geq e_3, \ \phi(s_2) > \phi(s_2-c\tau)$. 
However, in view of the inequalities 
$$
g(\phi(s_2),\phi(s_2-c\tau)) < g(\phi(s_2),\phi(s_2))  \leq 0,\  \mbox{if} \ \phi(s_2-c\tau) \geq \kappa; 
$$
$$
g(\phi(s_2),\phi(s_2-c\tau)) < g(\phi(s_2),e_1)  < 0,\  \mbox{if} \ e_1 < \phi(s_2-c\tau) < \kappa;
$$
$$
g(\phi(s_2),\phi(s_2-c\tau)) = g(\phi(s_2),e_1)  < 0,\  \mbox{if} \ e_1 = \phi(s_2-c\tau), 
$$
this contradicts to equation (\ref{twe2ng}). 

\noindent (b) Suppose that $\psi(s'-c\tau) \leq \kappa$ and $\psi'(s') =0$ for some $s' \leq s$. Then $\psi'(s')=0, \psi''(s')\leq 0$  so that 
$$0\leq  g(\psi(s'),\psi(s'-c\tau)) < g(\psi(s'-c\tau),\psi(s'-c\tau)) \leq 0, \  \mbox{if} \ \psi(s') \leq e_2  \  \mbox{and} \ \psi(s'-c\tau) <\psi(s'); $$
$$0\leq  g(\psi(s'),\psi(s'-c\tau)) = g(\psi(s'-c\tau),\psi(s'-c\tau)) < 0, \  \mbox{if} \ \psi(s') =\psi(s'-c\tau) \leq \kappa; $$
$$0\leq  g(\psi(s'),\psi(s'-c\tau)) < g(\psi(s'),e_1) <0, \  \mbox{if} \ \psi(s') > e_2, $$
a 
contradiction. The same argument works if we suppose that $\psi'(s')\geq 0, \psi''(s')=0$. \hfill \ter
\begin{lem}
\label{2af+}    Assume  {\bf (U)}  and suppose that $\phi(t) \geq e_1, \ t \in \R$,   is a bistable  wavefront for equation (\ref{twe2ng}) propagating with the speed $c>0$.  Then, for some appropriate $t_1\in \R$ and small $\epsilon >0$, 
$$
\phi(t+t_1) = e_3-e^{\mu_2 t} + O(e^{(\mu_2-\epsilon) t}),  \quad \phi'(t+t_1) = -\mu_2e^{\mu_2 t} + O(e^{(\mu_2-\epsilon) t}),  \ t \to +\infty. 
$$
In particular, $\phi(t)$ is eventually strictly increasing at $+\infty$. 
\end{lem}
{\it Proof. }  Recall that system (\ref{asy}) is  exponentially dichotomic  at $+\infty$. As a consequence, $\phi'(t), \phi''(t)$ and positive function $\phi_1(t)= e_3- \phi(t)=\int_t^{+\infty}\phi'(s)ds$
converge to $0$ exponentially at $+\infty$.  Thus, for some positive $\nu$ and $t\to +\infty$, functions  
$$
a_j(t) : = \int_0^1g_j(e_3-(1-s)\phi_1(t), e_3-(1-s)\phi_1(t-c\tau))ds 
$$
satisfy 
$$
a_1(t) = a_++ O(e^{-\gamma \nu t}),  \quad a_2(t) =b_+ + O(e^{-\gamma \nu t}), \ t \to +\infty.  
$$
Now, since positive $\phi_1(t)$ satisfies the equation
\begin{equation} \label{asyma*}
\phi_1''(t)-c\phi_1'(t)+a_1(t)\phi_1(t)+a_2(t)\phi_1(t-c\tau)=0, 
\end{equation}
the super-exponential convergence of $\phi_1(t)$ to $0$ as $t \to +\infty$ is not possible due  to  \cite[Lemma 3.1.1 under Assumption 3.1.1]{hl2}.  Thus Proposition 7.2 from \cite{MP} implies  that, for some eigenvalue $\mu_j, \ \Re\mu_j <0,$ of $\chi_+(z)$, small positive $r$ and non-zero polynomials $p_j(t), q_j(t)$, it holds 
$$
(\phi_1(t), \phi_1'(t)) = (p_j(t),q_j(t))e^{\mu_j t} + O(e^{(\Re \mu_j -r) t} ), \ t \to +\infty.  
$$
 In view of the positivity of $\phi_1$ at 
$+\infty$, $\mu_j$ should be  a real negative number. Thus actually 
$\mu_j = \mu_2$.  By Lemma \ref{lc2}, $\mu_2$ is a simple zero of $\chi_+(z)$ and therefore $p_j >0, q_j= \mu_2 p_j$ are constants.  \hfill \ter

In the next two lemmas, we study the asymptotic behavior of bistable wavefronts at $-\infty$. 
\begin{lem} \label{2af}  Assume   {\bf (U)} and suppose that  $\phi(t)$ is profile of a bistable wavefront which is monotone at  $-\infty$. Then, for some appropriate $t_1\in \R$ and small $\epsilon >0$, 
\begin{equation} \label{as22c} \hspace{-20mm}
\phi(t+t_1) = e_1+e^{\lambda_1 t} + O(e^{(\lambda_1+\epsilon) t}), \quad \phi'(t+t_1) = \lambda_1e^{\lambda_1 t} + O(e^{(\lambda_1+\epsilon) t}), \ t \to -\infty. 
\end{equation}
In particular,  $\phi'(t)>0$ on some maximal open interval $(-\infty,s)$, $s >0$,  described in Lemma \ref{cot}. 
\end{lem}
{\it Proof. }  In view of Lemma \ref{cot}, $\psi(t) = \phi'(t) > 0,\psi'(t)= \phi''(t) >0$  on some maximal interval $(-\infty,s)$. Since  ${\mathcal D}\psi(t)=0$, we find  that 
$(\psi'(t)e^{-ct})' >0$ because of 
$$
\psi''(t)- c \psi'(t)  =-a(t)\psi(t) - b(t)\psi(t-c\tau) >0, \quad t < -c\tau.  
$$
Therefore $\psi'(t)e^{-ct} < \psi'(s)e^{-cs}$ for $t<s <-c\tau$, or, equivalently, $\phi''(t) < \phi''(s)e^{-c(s-t)}$ for $t<s<-c\tau$. 
Thus $\phi'(t) = \int_{-\infty}^t\phi''(s)ds$, 
$\phi''(t)$ converge to $0$ exponentially at $-\infty$, so that, for some positive $\nu$,  
$$
a(t) = a_-+ O(e^{\nu t}),  \quad b(t) = b_-+ O(e^{\nu t}), \quad t \to - \infty. 
$$
Applying now \cite[Lemma 3.1.1 under Assumption 3.1.2]{hl2},  Proposition 7.2 from \cite{MP}, we obtain that, for some eigenvalue $\lambda_j, \ \Re\lambda_j >0,$ of $\chi_-(z)$, small positive $r$ and non-zero polynomials $p_j(t), q_j(t)$, it holds 
$$
(\phi'(t), \phi''(t)) = (p_j(t),q_j(t))e^{\lambda_j t} + O(e^{(\Re \lambda_j +r) t} ), \ t \to -\infty.  
$$
Now, eventual  monotonicity of $\phi(t)$ at $-\infty$ implies eventual non-negativity or non-positivity of $\phi'(t)$. 
Thus  $\lambda_j$ should be  a real positive number.  This yields that actually 
$\lambda_j = \lambda_1$ and $p_j$ is a non-zero constant. Now, if $p_j$ is negative, then $\phi(t)$ is strictly decreasing at $-\infty$ and therefore there is the leftmost number $s$ such that $\phi'(s) =0$, $\phi''(s) \geq 0,$ 
$e_1 > \phi(s-c\tau)> \phi(s)$. This contradicts, however, to equation (\ref{twe2ng}): 
$$
0 \geq g(\phi(s),\phi(s-c\tau)) > g(\phi(s),e_1) >0.
$$
 In this way,  $p_j >0$ that  proves the second formula in (\ref{as22c}) for an appropriate $t_1$, while the similar formula for $\phi(t)$ at $t\to -\infty$ \ follows from the representation $\phi(t) - e_1= \int^t_{-\infty}\phi'(s)ds$. 
\hfill \ter

The next  result (used immediately afterwards, in the proof of Lemma \ref{2aaf}) excludes the existence of small solutions to asymptotically autonomous delayed differential equations at $-\infty$: 
\begin{lem} \label{sms} Suppose that $L,M(t):C([-h,0],\R^n) \to \R^n,\ t \leq 0$, are continuous linear operators and $\|M(t)\| \to 0$ as $t\to -\infty$
(here $\|\cdot\|$ denotes the operator norm).  Then the system 
\begin{equation} \label{LM}
x'(t) = (L+M(t))x_t, \  x_t(s):= x(t+s), \ s \in [-h,0],
\end{equation}
does not have exponentially small solutions at $-\infty$ (i.e. non-zero solutions $x:\R_- \to \R^n$ such that for each $\gamma \in \R$ it holds that $x(t)e^{\gamma t} \to 0, \ t \to -\infty$). 
\end{lem}
{\it Proof. }  On the contrary, suppose that there exists a small solution $x(t)$ of $(\ref{LM})$ at $-\infty$. Take some $b>0$.  It is straightforward  to see that the property 
 $x(t)e^{\gamma t} \to 0,$  $t \to -\infty,$ is equivalent to  $|x_t|_be^{\gamma t} \to 0, \ t \to -\infty$, where $|x_t|_b = \max_{s \in [-b,0]}|x(t+s)|$. Next, smallness of 
 $x(t)$ implies that 
 $
 \inf_{t \leq 0} |x_{t-b}|_b/|x_t|_b =0. 
 $
 Indeed, otherwise there is $K>0$ such that $|x_{t-b}|_b/|x_t|_b \geq K, \ t \leq 0,$ and therefore, setting $\nu:= b^{-1}\ln K$, we obtain the following contradiction: 
 $$
0<  |x_t|_be^{\nu t} \leq |x_{t-b}|_be^{\nu (t-b)}\leq |x_{t-2b}|_be^{\nu (t-2b)}\leq\dots\leq |x_{t-mb}|_be^{\nu (t-mb)} \to 0, \quad m \to +\infty.  
$$
Hence, for $b=3h$  there is  a sequence $t_j\to -\infty$ such that  $|x_{t_j-3h}|_b/|x_{t_j}|_{3h} \to 0$ as $j \to \infty$. Clearly, 
$|x_{t_j}|_{3h} = |x(s_j)|$ for some $s_j \in [t_j-3h,t_j]$ and, for all large $j$, it holds $|x(s_j)| \geq |x(s)|, \ s \in [t_j-6h,t_j]$.   Since $0\leq t_j-s_j \leq 3h$, 
without loss of generality we can assume that $\theta_j:= t_j-s_j \to \theta_* \in [0,3h]$. 

Now, for sufficiently large $j$, consider the sequence 
of functions $$
y_j(t)= \frac{x(t+t_j)}{|x(s_j)|}, \ t \in [-6h,0], \quad |y_j(-\theta_j)| =1, \quad  |y_j(t)| \leq 1, \ t \in [-6h,0]. 
$$
For each $j$, $
y_j(t)$ satisfies the equations 
$$
y'(t) = (L+M(t+s_j))y_t, \quad y_j(t) = y_j(-\theta_j)+ \int^t_{-\theta_j}(L+M(u+s_j))y_udu, 
$$
and therefore   
$
 |y_j(t)| \leq 1, \ |y'(t)| \leq \|L\|+\sup_{s \leq 0} \|M(s)\|, \ t\in  [-5h,0],\ j \in \N.
$
Thus, due to the Arzel\`a-Ascoli theorem,  there exists a subsequence $y_{j_k}(t)$ converging, uniformly on  $[-5h,0]$, to some continuous function $y_*(t)$ such that 
$|y_*(-\theta_*)|=1, $
$$
y_*(t) = y_*(-\theta_*)+ \int^t_{-\theta_*}L(y_*)_udu,\quad  t\in  [-4h,0]. 
$$
In particular, $y'_*(t) = L(y_*)_t,\  t\in  [-4h,0]$.  Since $y_*(t)= 0$ for all $t\in [-5h,-3h]$,  the existence and uniqueness theorem applied to the initial 
value problem $y'(t) = Ly_t,\  t\in  [-4h,0], \ y_{-3h} =0,$ implies that  also $y_*(t)= 0$  for all $t \in [-3h,0]$. However, this contradicts  the relation
$|y_*(-\theta)|=1$. The proof of Lemma \ref{sms} is completed. 
\hfill \ter
\begin{lem} \label{2aaf} Assume condition {\bf (U)}.   If $\phi(t)$ is a  solution of equation (\ref{twe2ng}) such that $\phi(0) >e_1$ and 
\begin{equation}\label{lam}
\sup_{t \leq 0}|\phi(t)-e_1|e^{-(\lambda_1-\delta)t} < \infty
\end{equation}
for some $\delta \in (0, \lambda_1)$ small enough to satisfy 
$\Re \lambda_2 < \lambda_1-1.5\delta$, 
 then $\phi'(t)>0$ on the maximal open interval $(-\infty,s)$ described in Lemma \ref{cot}. 
\end{lem}
{\it Proof. }  Set $\phi_1(t)= e_1-\phi(t)$ and $\nu = \lambda_1-\delta$. Then  $\phi_1(t)$ satisfies  equation (\ref{asyma*}) where
$$
a_j(t) : = \int_0^1g_j(e_1-(1-s)\phi_1(t), e_1-(1-s)\phi_1(t-c\tau))ds, 
$$
$$
a_1(t) = a_-+ O(e^{\gamma \nu t}),  \quad a_2(t) =b_- + O(e^{\gamma \nu t}), \ t \to -\infty.  
$$
By Lemma \ref{sms}, $\phi_1(t)$  is not super-exponentially small  at $-\infty$.  Since $\phi_1'(-\infty) =0$, we can again invoke Proposition 7.2 from \cite{MP} and Lemma \ref{lc3}  to conclude that, for some eigenvalue $\lambda_j, \ \Re\lambda_j >0,$ of $\chi_-(z)$, small positive $r$ and non-zero constants  $p_j, q_j$, it holds that 
$$ 
(\phi_1(t), \phi_1'(t)) = (p_j,q_j)e^{\lambda_j t} + O(e^{(\Re \lambda_j +r) t} ), \ t \to -\infty.  
$$
 But then, in the latter case,  condition (\ref{lam}) implies that $\lambda_j = \lambda_1$.  
This leads to the conclusion of the lemma, see  the final part of the proof of Lemma \ref{2af}.  
\hfill \ter

In the remainder of this section, we are assuming that the profile $\phi(t)$ of bistable wave is strictly increasing, i.e.  
$\phi'(t) >0, \ t \in \R$.    
After fixing some $\delta \in (0, 0.5(\lambda_1-\Re \lambda_2))\cap (0, \lambda_1)\cap (0, -\mu_2)$, such that 
$(1+\gamma)(\mu_2+\delta) < \mu_2$, 
we will consider ${\frak F}_c$ as a linear operator defined
on $C^1_\delta$ and taking its values in $C_\delta$, where 
\[\hspace{-25mm}
C_\delta = \{y= (y_1,y_2) \in C(\R, \R^2): |y|_\delta: = \sup_{s \leq
0}e^{-(\lambda_1- \delta)s} |y(s)|+  \sup_{s \geq
0}e^{-(\mu_{2} + \delta)s} |y(s)| < \infty  \},
\]
\[
C^1_\delta = \{y \in C_\delta: y, y' \in C_\delta, \
|y|_{1,\delta}: =  |y|_\delta +  |y' |_\delta < +\infty
\}. \]
\begin{remark}\label{pm} From Lemmas \ref{2af+} and \ref{2af} we obtain  that $(\phi',\phi'') \in C_\delta$ and that
$
a(t) = a_++ O(e^{\gamma \mu_2 t}), $ $b(t) = b_++ O(e^{\gamma\mu_2 t}), \quad t \to +\infty. 
$
\end{remark} 
The system which is formally adjoint \cite{hale} to (\ref{asy}) has the following  form 
\begin{equation}\label{asy*}
  v'(t) = a(t)w(t)+ b(t+c\tau)w(t+c\tau), \ \
   w'(t)  = -v(t)-cw(t).
\end{equation}
This amounts to the following equation for $w(t)$:
$$
{\mathcal D}^*w(t): =w''(t)+cw'(t) + a(t)w(t) + b(t+c\tau)w(t+c\tau) =0.
$$
The following result is obvious. 
\begin{lem} \label{nt} Suppose that functions $a(t), b(t)$ are continuous and $b(t)\not= 0$ if $t\not= 0$. If ${\mathcal D}^*w(t) =0$ and, for some $t'$, it holds that either $w(t)=0$ for all $t \leq t'$ or $w(t)=0$ for all $t \geq t'$,  then 
$w(t)\equiv 0$. 
\end{lem}
The equation ${\mathcal D}^*w(t) =0$ with advanced argument can be transformed in the usual delayed equation by means of  the transformation $w(t) \to  w(-t)$:  
\begin{equation}\label{rev}
w''(t)-cw'(t) + a(-t)w(t) + b(-t+c\tau)w(t-c\tau) =0.
\end{equation}
A  description of the leading zeros of  the  characteristic functions
$\chi_\pm(z) $ provided in Lemmas \ref{lc2}, \ref{lc3} together with  an analysis of 
the asymptotic properties of wavefronts realized in Lemmas \ref{cot}, \ref{2af+}, \ref{2af}, \ref{2aaf} and Remark \ref{pm} enable us 
to establish the Fredholm properties of operators ${\mathcal D}, {\frak F}_c$.  Note that Fredholmness of asymptotically autonomous  functional differential 
operators (respectively, of delayed, of mixed or of non-local type and considered in different spaces) was studied by Hale and Lin \cite{HL},  Mallet-Paret \cite{MP} and in \cite{ADV,DMV,V4}. In fact, the next result can be deduced from any of  these works: 
\begin{pro} \label{P1+} ${\frak F}_c:C^1_\delta  \to C_\delta$ is Fredholm operator of index ind ${\frak F}_c=0$. Moreover, ${\frak F}_c$ has one-dimensional  kernel $N({\frak F}_c) =< (\phi'(t), \phi''(t))>$. Thus range $R({\frak F}_c)$ of ${\frak F}_c$ has codimension 1 and therefore 
$$
R({\frak F}_c)= \{F=(f_1,f_2) \in C_\delta: \int_\R \left(f_1(s)v_*(s)+ f_2(s)w_*(s)\right)ds =0\}, 
$$
where $z(t)= (v_*(t),w_*(t))$ is the unique (up to a constant multiple) non-zero solution of  (\ref{asy*}) satisfying inequalities 
$$
|z(t)| \leq Ke^{-(\mu_1-0.5\delta)t}, \ t \geq 0; \quad |z(t)| \leq Ke^{-(\lambda_1-1.5\delta) t}, \ t \leq 0. 
$$
\end{pro}
{\it Proof. }  Due to our choice of spaces, it is convenient to use the theory developed in \cite{HL}.  Particularly,  we are going to show how 
Proposition \ref{P1+} can be deduced from  Lemmas 4.6 and 4.5 in \cite{HL}.  In order to simplify the use of these lemmas and to show their relation to similar results in \cite{ADV,DMV,MP},  for $(p,q)= (\lambda_1-\delta,\mu_2 + \delta)$ , $q\leq 0\leq p,$ we consider $C^\infty$-smooth weight function $\mu: \R\to [1,+\infty)$ such that $\mu(t) = e^{-pt}, \ t \leq -1$ and 
$\mu(t) = e^{-qt}, \ t \geq 1$. Following \cite{HL}, we introduce the notation  
\[
C(p,q) = \{y= (y_1,y_2) \in C(\R, \R^2): |y|_{p,q}: = \sup_{s \in \R}\mu(s)|y(s)| < \infty  \},
\]
\[
C^1(p,q) = \{y \in C(p,q): y, y' \in C(p,q), \
|y|_{1,p,q}: =  |y|_{p,q} +  |y' |_{p,q} < +\infty
\}. \]
Clearly, $C_\delta = C(\lambda_1-\delta,\mu_2 + \delta)$, $C^1_\delta = C^1(\lambda_1-\delta,\mu_2 + \delta)$ and the multiplication operator $Ly = \mu y$, 
$$ L: C^1(p,q) \to  C^1(0,0), \quad L:  C(p,q) \to  C(0,0), $$
 is an isomorphism of the Banach spaces. Consider $B: = L{\frak F}_c L^{-1}:  C^1(0,0) \to  C(0,0)$, 
 \[ B(v,w)(t)
= (v'(t)-\frac{\mu'(t)}{\mu(t)}v(t)-w(t), \]\[w'(t)-\frac{\mu'(t)}{\mu(t)}w(t)+  a(t)v(t)- c w(t) + b(t)\frac{\mu(t)}{\mu(t-c\tau)}v(t-c\tau)).
\]
Since the limiting equations $B_{\pm\infty}(v,w)(t) =0$ at $-\infty$ and $+\infty$, 
$$
B_{-}(v,w)(t)
= (v'(t)+pv(t)-w(t), w'(t)+pw(t)+  a_-v(t)- c w(t) + b_-e^{-pc\tau}v(t-c\tau)).
$$
$$
B_{+}(v,w)(t)
= (v'(t)+qv(t)-w(t), w'(t)+qw(t)+  a_-v(t)- c w(t) + b_-e^{-qc\tau}v(t-c\tau)),
$$
have the characteristic functions $\chi_-(z+p)= \chi_-(z+\lambda_1-\delta), \ \chi_+(z+q)= \chi_+(z+\mu_2+\delta)$, 
they both are exponentially dichotomic with one-dimensional unstable spaces, cf. Lemma \ref{lc2},\ref{lc3}.  Then 
Lemmas 4.6 and 4.5 in \cite{HL} imply that $B$ (and, consequently, ${\frak F}_c= L^{-1}BL$) is Fredholm of index $0$. The same conclusion can be obtained by using approach proposed in  \cite[Theorems 2.2 and 3.2]{DMV}. Next, it is clear (see also \cite[Lemma 4.6]{HL}) that dimension of kernel $N(B) $ of $B$ cannot exceed the dimension of the unstable space of $B_-$, i.e $\dim N(B) \leq 1$. 
On the other hand, by Remark \ref{pm}, we have that $L(\phi',\phi'')=\mu(\phi',\phi'') \in C(0,0)\cap N(B) $ so that $\dim N(B) = \dim N({\frak F}_c)={\rm codim}\, R({\frak F}_c)= {\rm codim}\, R(B) =1$. 

Finally,  \cite[Lemma 4.5]{HL} (or, equivalently,  \cite[Theorem A and Proposition 5.3]{MP}) implies that 
$R(B)= \{F \in C(0,0): \int_\R F(s)y_*(s)ds =0\}$, where $y_*(s)$ is an exponentially decaying non-zero solution of the formally adjoint  equation $B^*y=0$, 
 \[ B^*(v,w)(t)
= (v'(t)+\frac{\mu'(t)}{\mu(t)}v(t)-a(t)w(t)- b(t+c\tau)\frac{\mu(t+c\tau)}{\mu(t)}w(t+c\tau), \]\[w'(t)+\frac{\mu'(t)}{\mu(t)}w(t)+  v(t) +c w(t)).
\]
The  characteristic functions of the limiting equations for $B^*$ at $-\infty$ and $+\infty$ are, respectively, $\chi_-(p-z)= \chi_-(\lambda_1-\delta-z), \ \chi_+(q-z)= \chi_+(\mu_2+\delta-z)$. Therefore,  
from  \cite[Proposition 7.2 ]{MP} we find that $y_*(t) = O(te^{(\lambda_1-\lambda_2 -\delta)t})$ at $-\infty,$ and  $y_*(t) = O(e^{(\mu_2-\mu_1 +\delta)t}), \ t \to +\infty.$
Thus $$R({\frak F}_c) = L^{-1}R(B) = \{F \in C_\delta: \int_\R F(s)\mu(s)y_*(s)ds =0\} =  \{F \in C_\delta: \int_\R F(s)z(s)ds =0\},$$
where $z(s):= \mu(s)y_*(s)$ clearly satisfies the estimates of the proposition as well as equation (\ref{asy*}); the latter can be verified by a direct calculation. 
\hfill \ter
\begin{remark} \label{R1} In view of Lemma \ref{nt}, we can assume that $w_*(s)>0$ for some $s >0$. Next, set $w(t)= w_*(-t)$. Then  $w'(t)= -w_*'(-t) =  v_*(-t)+cw_*(-t)$ so that 
$$
|(w(t),w'(t))| = O(e^{(\mu_1-0.5\delta)t}), \ t \leq 0; \quad |(w(t),w'(t))| = O(e^{(\lambda_1-1.5\delta) t}), \ t \geq 0. 
$$
Since the characteristic function of the limiting equation for (\ref{rev}) at $+\infty$ is $\chi_-(z)$, Proposition   7.2 in \cite{MP} guarantees that either $w(t)$ is  super-exponentially small at $+\infty$ or, for some eigenvalue 
$\lambda_j$ with the real part $\Re \lambda_j \leq \Re \lambda_2$ and for some polynomial $P(t)$ of the degree less than or equal to one,  it holds 
\begin{equation}\label{str}
w(t) = P(t) e^{\Re \lambda_j t}\cos (\Im \lambda_j t +\phi) + O( e^{(\Re \lambda_j-\varepsilon) t}), t \to +\infty.  
\end{equation}
{In particular, this shows that either $w_*(t)$ is exponentially (or even super-exponentially) small at $-\infty$ or $w_*(t)$ is non-decaying and oscillating around $0$ at $-\infty$. }
\end{remark}
\begin{cor} \label{Co2} Set $Y_\delta: = \{y \in  C(\R, \R):  (y,y) \in C_\delta\}$ and  $X_\delta: = \{y \in  C^2(\R, \R):  y, y',y'' \in Y_\delta\}$, then ${\mathcal D}:X_\delta \to Y_\delta$ is continuous Fredholm operator of index 0 and one-dimensional kernel $N({\mathcal D}) =< \phi'(t)>$.  The range $R({\mathcal D})$ of ${\mathcal D}$ is given by 
$$
R({\mathcal D})= \{f \in Y_\delta: \int_\R f(s)w_*(s)ds =0\},
$$
where $w_*(t)$ is described in Proposition \ref{P1+}. 
\end{cor}
{\it Proof. }  Indeed, $y \in N({\mathcal D})$ if and only if $(y,y') \in N({\frak F}_c)$. Similarly, $f \in R({\mathcal D})$
if and only if $(0,f) \in R({\frak F}_c)$, i.e. if and only if  $\int_\R f(s)w_*(s)ds =0$. 
\hfill \ter
\begin{lem} \label{Lp}  The solution $w_*$ is positive on $(0,+\infty)$:  $w_*(t) > 0, \ t > 0$. 
\end{lem} 
{\it Proof. }  To prove the non-negativity of $w_*(t)$ on $\R_+$,  we are going to use, similarly to the proofs of  \cite[Theorem 5.1]{VVV} or \cite[Theorem 2.5]{hl1},  an appropriate test function $f$. 
 Recall that, by our assumption (see Remark \ref{R1}),  $w_*(s)>0$ for some $s >0$. 

\vspace{2mm}
 
\noindent \underline{Claim 1.} {\it  It holds that  $w_*(t) \geq 0$ for all $t\in \R_+$.}\\
 Indeed, otherwise we can indicate 
a function $f \in Y_\delta$ and a real number $T>0$ with the following properties
\begin{enumerate}
\item[(i)]  $f(t)=0, \ t \leq 0,$ and $f(t) <0, \ t >0; $
\item[(ii)]  $f(t) = {\mathcal D}(te^{\mu_2 t}) = e^{\mu_2 t}\left(\chi_+'(\mu_2)  +o(1)\right) <0, \ t \geq T;$
\item[(iii)]  $f(t)$ is smooth on $\R_+$ and $f'(0+) <0$;
\item[(iv)]  $\int_\R f(t)w_*(t)dt =0$. 
\end{enumerate}
Then, by Corollary \ref{Co2},  the inhomogeneous equation 
\begin{equation}\label{test*}
\psi''(t)-c \psi'(t) + a(t)\psi(t) + b(t)\psi(t-c\tau) = f(t)
\end{equation}
has a solution $\psi_* \in X_\delta$.  
Since $f(t) =0$ for $t\leq 0$, we conclude that there exists $T_2>0$ such that the vector  $(\psi_*(t),\psi_*'(t)), \ t \leq -T_2, $ belong to the unstable space of the  system ${\frak F}_c(v,w)=0, \ t \leq -T_2$ which has a shifted exponential dichotomy at $-\infty$.  Now, since  the exponents $\alpha_1:=\lambda_1-1.5\delta < \lambda_1-0.5\delta=: \beta_1 < \lambda_2$ of the shifted exponential dichotomy  satisfy
$\Re \lambda_2 < \alpha_1< \beta_1 < \lambda_1$,  
this unstable space has dimension 1 and therefore  $\psi_*(t) = k\phi'(t), \ t \leq -T_2,$ for some $k>0$. For certain, this yields immediately that $\psi_*(t) = k\phi'(t), \ t \leq 0$. 
 
 On the other hand,  due to our definition of  $f(t)$ for positive $t$, we obtain that the function $q(t):=\psi_*(t)- te^{\mu_2 t}, \ t \geq  T + c\tau$ satisfies the homogeneous equation  ${\mathcal D}q(t)=0, \ t \geq  T + c\tau$. 
Since $q(t), q'(t)$ have an exponential rate of convergence to $0$ at $+\infty$, we can apply Proposition 7.2 from \cite{MP} together with Remark \ref{pm}, to conclude that $q(t) = O(e^{\mu_2t}), \ t \to +\infty$.  This shows that 
the solution 
$
\psi(t,\xi) = \psi_*(t) + \xi \phi'(t)
$
of (\ref{test*}) is positive at $+\infty$ for every real $\xi$.  Since also $\psi(t,\xi) = (k+\xi)\phi'(t), \ t \leq 0$, we obtain that 
$
\psi(t,\xi) >0, \  t \in \R, 
$
for all large $\xi >0$.  Let now $\xi_*: = \inf\{\xi: \psi(t,\xi) >0, \  t \in \R\}$. Clearly, $\xi_*$ is finite and 
$\psi(t,\xi) \geq 0, \ t \in \R,$ if and only if $\xi \geq \xi_*$.  Next, we have that either $\psi(t,\xi_*)>0$ for $t\leq 0$ 
or $\psi(t,\xi_*)\equiv 0$ on $\R_-$. In the first case, $\psi(t_*,\xi_*)=0$ for some $t_* >0$ (since otherwise 
$\psi(t,\xi_*-\epsilon)>0, \ t \in \R$, for all small $\epsilon >0$).  However, this  implies that 
$\psi''(t_*,\xi_*)\geq \psi'(t_*,\xi_*)=0$  contradicting to (\ref{test*}) at $t_*$ (since $f(t_*)<0$, $b(t_*) >0$ and $\psi(t_*-c\tau,\xi_*)\geq 0$). In the second case, $C^2$-smooth function $\psi(t,\xi_*)$ satisfies on $[0,c\tau]$ the following ordinary equation with zero initial data: 
$$
\psi''(t)-c \psi'(t) + a(t)\psi(t) = f(t), \ \psi(0)= \psi'(0) =0.
$$
In particular,  $\psi''(0,\xi_*)=0$. However, $\psi'''(0+,\xi_*)= f'(0+)<0$ and therefore $\psi(t,\xi_*)<0$ for all small 
positive $t$, contradicting to the definition of $\xi_*$.  This completes the proof of Claim 1.

\vspace{2mm}

\noindent
\underline{Claim 2.} {\it  It holds that  $w_*(t) > 0$ for all $t >0$.}\\
To analyse the asymptotical behaviour of $w_*(t)$ for positive $t$, it is convenient to consider $\hat w(t)=w_*(-t)$. This function 
satisfies the delayed equation (\ref{rev})  which is asymptotically autonomous at $-\infty$, with the limiting equation 
$$
w''(t)-cw'(t) + a_+w(t) + b_+w(t-c\tau) =0.
$$
Since $|\hat w(t)| \leq Ke^{(\mu_1-0.5\delta)t}, \ t \leq 0$, and by Remark \ref{pm} $a(-t)-a_+,\ b(-t+c\tau)-b_+$ are exponentially small at $-\infty$,  we deduce,  as before,  from  \cite[Proposition   7.2]{MP} and  Lemma \ref{sms}  the following asymptotic representation 
$$
\hat w(t) = d e^{\mu_1 t} + O (e^{(\mu_1+\epsilon)t}),   \quad t \to -\infty,
$$
with some positive $\epsilon, d$. Consequently,  $w_*(t) >0$ on some interval $(m, +\infty)$. 
Let $m$ be the leftmost point for which the inequality  $w_*(t) >0,\ t \in (m, +\infty)$ holds. 
If $m >0$, then  $w_*(m)=w_*'(m)= 0 \leq w_*''(m)$. Therefore, in view of equation 
 ${\mathcal D}^*w_*(t)=0$, we find that $b(m+c\tau)w(m+c\tau) \leq 0$, contradicting to the fact 
 that $b(m+c\tau) >0, \ w(m+c\tau) >0$. 

\hfill \ter
\begin{remark} The second example considered in Subsection  \ref{toym}  (with $k_*>1$ and $c <0$) shows that, in general, $w_*(t)$ can oscillates on $(-\infty,0)$.  
On the other hand, we believe that $w_*(t) >0, \ t \in \R$, if  $|c(\tau)| < clin(\tau)$, cf. Fig.3.  See the next section where we prove such a kind of result
when {\bf (U$^*$)} is assumed instead of {\bf (U)}. 
\end{remark} 

\section{Variational equation along the monotone bistable wave, case of hypotheses {\bf (B),  (U$^*$)}.}\label{S4}
Let profile $\phi(t)$ of the bistable wavefront for problem (\ref{twe2ng}) considered with $c>0$ be such that   $\phi'(t) >0, \ t \in \R$.   Again,  we can assume that $\phi(-c\tau)= \kappa$ and that $\phi(t) < \kappa$ for $t < -c\tau$.  In view of assumptions   {\bf  (B),  (U$^*$)},  the coefficients $a(t), \ b(t)$ of the differential operator ${\mathcal D}$ satisfy the relations 
$
b(0)=0,$ \ $ b(t) >0$ for  $t < 0$  and $ a(t), \ b(t) < 0$ for  $t > 0$, while 
$$
a_+:=a(-\infty)= g_1(e_1, e_1) <0, \quad b_+:=b(-\infty) = g_2(e_1, e_1) > 0;  \quad b_++a_+ <0; 
$$
$$
 a_-:=a(+\infty)= g_1(e_3, e_3) { < 0}, \quad b_-:= b(+\infty) = g_2(e_3, e_3) <0.
$$
Hence, the variational equation is asymptotically autonomous and the limiting autonomous equations at $\pm\infty$
have the  characteristic functions
$
\chi_\mp(z). 
$
Clearly, above convention on the notation allows the application of  Lemmas \ref{lc2} and \ref{lc3} describing properties of zeros of $\chi_\mp(z)$.  
In Lemma \ref{2aff} below, we show how the monotonicity of  wavefront $\phi(t)$ at $+\infty$ propagating with speed $c$  implies that $\chi_-(z)$ has exactly two (counting the multiplicity) real negative zeros $\lambda_3 \leq \lambda_2 <0$ (i.e. implying that $(\tau,c) \in \frak{ D}$). Therefore $\phi'(t)$ decays 
at $+\infty$ with the exponential rate which is asymptotically equivalent  to $p(t)\exp(\lambda_jt)$, where $j \in \{2,3\}$ and $p(t)$ is a polynomial.  Our approach, however, requires slowest
possible decay of $\phi'(t)$ at $+\infty$. We are reaching this goal assuming the sub-tangency 
condition at the steady state $e_3$ in the hypothesis ({\bf U$^*$}). As we show in  Lemma \ref{2aff}, this condition forces  $\phi'(t)$ to have the required asymptotical behavior  at $+\infty$. It is worth to mention that the slowest
decay of $\phi'(t)$ at $+\infty$ was automatically assured in the case of  the hypothesis ({\bf U}). This explains why 
a similar sub-tangency condition was not required in ({\bf U}). 
\begin{lem} \label{2aff}   Let the hypotheses {\bf (B),  (U$^*$)} be satisfied. If equation  (\ref{twe2an}) has a non-decreasing 
bistable wavefront and  $\chi_-(z)$ does not have roots on the imaginary axis, then $\chi_-(z)$ has exactly two negative zeros (counting multiplicity) $\lambda_3\leq \lambda_2$.  Moreover,  for some appropriate $t_1\in \R$, $A >0, \ j \in \{0,1\},$  and small $\epsilon >0$, 
\begin{equation} \label{as22}\hspace{-25mm}
\phi(t+t_1) = e_3 - At^je^{\lambda_2 t } + o(t^je^{\lambda_2 t} ), \ \phi'(t+t_1) = - A\lambda_2t^{j}e^{\lambda_2 t} + o(t^je^{\lambda_2 t}), \ t \to +\infty. 
\end{equation}
Here $j=1$ if and only if $\lambda_2= \lambda_3$. 
\end{lem}
{\it Proof. } Since  $\chi_-(z)$ does not have roots on the imaginary axis, system (\ref{asy}) is exponentially dichotomic  at $+\infty$. As a consequence, $\phi'(t), \phi''(t)$ converge to $0$ exponentially fast at $+\infty$.  Thus, for some positive $\nu$,  
$$
a(t) = a_-+ O(e^{-\nu t}),  \quad b(t) = b_-+ O(e^{-\nu t}), \quad t \to +\infty. 
$$
Applying now \cite[Lemma 3.1.1 under Assumption 3.1.2]{hl2},  Proposition 7.2 from \cite{MP}, we obtain that, for some eigenvalue $\lambda_j, \ \Re\lambda_j < 0,$ of $\chi_-(z)$, small positive $r$ and non-zero polynomials $p_j(t), q_j(t)$, it holds that
\begin{equation} \label{fn*}
(\phi'(t), \phi''(t)) = (p_j(t),q_j(t))e^{\lambda_j t} + O(e^{(\Re \lambda_j -r) t} ), \ t \to +\infty.  
\end{equation}
Now, monotonicity of $\phi(t)$ at $+\infty$ implies non-negativity  of $\phi'(t)$. 
Thus  $\lambda_j$ should be  a real negative number.  This yields that actually 
$\lambda_j \in \{\lambda_2, \lambda_3\}$ and $p_j$ is a positive constant (if $\lambda_3 < \lambda_2$) or at most first order non-zero polynomial 
(if $\lambda_2 = \lambda_3$).  The similar formula for $\phi(t)$ at $t\to +\infty$ \ follows from the equality $\phi(t) - e_3= -\int_t^{+\infty}\phi'(s)ds$. 

In order to prove that $j=2$ in the case when $g$ is sub-tangential at $e_3$, we observe that function 
$y(t):= \phi(t)-e_3$ satisfies the equation 
\begin{equation} \label{WWa}
y''(t) -cy'(t) + a_-y(t)+b_-y(t-c\tau) = h(t), \quad t \in \R,
\end{equation}
where 
$h(t) = a_-(\phi(t)-e_3)+b_-(\phi(t-c\tau)-e_3) - g(\phi(t),\phi(t-c\tau)) \geq 0$ because of the assumed sub-tangency of $g$. Furthermore, $h(t)\not\equiv 0$ since otherwise $y(t)$ should be equal to $0$, as a unique bounded solution of the exponentially dichotomic equation.   Since $g \in C^{1,\gamma}$, we also obtain that $h(t)= O(t^{1+\gamma} e^{(1+\gamma)\lambda_j t})$ at $t=+\infty$.  Therefore, applying the bilateral Laplace transform approach to equation (\ref{WWa}), we find that, for some small positive $r>0$, it holds
$$
y(t)= \mbox{Res}_{z=\lambda_2} \frac{e^{zt}\tilde h(z)}{\chi_-(z)} + O(e^{(\lambda_2-r) t}), \quad t \to +\infty.
$$
Here $\tilde h(z) = \int_\R e^{-zs}h(s)ds,\  \Re z \in ((1+\gamma)\lambda_2, 0), $ is the bilateral Laplace transform of $h(t)$. 
A simple calculation shows that 
$$
\mbox{Res}_{z=\lambda_2} \frac{e^{zt}\tilde h(z)}{\chi_-(z)}= \frac{e^{\lambda_2t}\tilde h(\lambda_2)}{\chi'_-(\lambda_2)} = -Ae^{\lambda_2t}, \ A:=  -\frac{ \int_\R e^{-\lambda_2 s}h(s)ds}{\chi'_-(\lambda_2)}>0, \quad \mbox{if \ } \lambda_3 < \lambda_2; 
$$
$$
\mbox{Res}_{z=\lambda_2} \frac{e^{zt}\tilde h(z)}{\chi_-(z)}=  -(Bt +D)e^{\lambda_2t}, \ B:=  -\frac{ 2\int_\R e^{-\lambda_2 s}h(s)ds}{\chi''_-(\lambda_2)}>0, \quad \mbox{if \ } \lambda_3 = \lambda_2.  
$$
This completes the proof of the lemma. 
\hfill \ter
\begin{remark} \label{W1} If we take $(\tau, c) \in \frak{ D}(\tilde a_-,\tilde b_-) \subset \frak{ D}(a_-,b_-) $ where $\tilde a_-,\tilde b_-$ were defined in (\ref{abt}), then 
the proof of Lemma \ref{2aff} works even without the sub-tangency condition. Indeed, $(\tau, c) \in \frak{ D}(\tilde a_-,\tilde b_-)$ implies that  $\chi_-(z)$ has exactly two negative zeros (counting multiplicity) $\lambda_3\leq \lambda_2$. Thus we obtain the following 
assertion. 
\begin{lem} \label{2affW}   Let the hypotheses {\bf (B),  (U$^*$)} (without the sub-tangency condition) be satisfied and  $(\tau, c) \in \frak{ D}(\tilde a_-,\tilde b_-)$.  If equation  (\ref{twe2an}) has a non-decreasing 
bistable wavefront, then    for some appropriate $t_1\in \R$, $A >0, \ j \in \{0,1\},$  and small $\epsilon >0$, the representation (\ref{as22}) is valid.  
\end{lem}
{\it Proof. }  First, suppose that, given  $(\tau, c) \in \frak{ D}(\tilde a_-,\tilde b_-)$, we have that  $\lambda_3 = \lambda_2$. Since this equality can occur only for $(\tau, c)$ on the 
boundary of domain $\frak{ D}(a_-,b_-)$, we conclude  that $(\tilde a_-,\tilde b_-) = (a_-,b_-)$.  Since this situation was already analyzed in Lemma \ref{2aff}, 
 we have to consider the case $\lambda_3 < \lambda_2$ only.  Consequently,  if the formula (\ref{as22})  does not hold, then 
it should be replaced with 
\begin{equation}\label{40}
\phi(t+t_1) = e_3 - e^{\lambda_3 t } + o(e^{\lambda_3 t} ), \quad \phi'(t+t_1) = -\lambda_3e^{\lambda_3 t} + o(e^{\lambda_3 t}), \ t \to +\infty. 
\end{equation}
Then replacing in equation (\ref{WWa}) $(a_-,b_-)$ with $(\tilde a_-,\tilde b_-)$ and arguing 
as below (\ref{WWa}),  we find that  $h(t)= O(e^{\lambda_3 t})$ at $t=+\infty$.  Applying the bilateral Laplace transform method again,  we obtain then that $\phi(t) =  e_3 - (At+B)e^{\tilde \lambda_2t}(1+o(1)),\  t \to+\infty,$ where $A,B$ satisfy $|A|+|B|>0$ and 
 $\tilde \lambda_2 \in (\lambda_3, \lambda_2)$ is the largest negative root of the equation $z^2-cz+\tilde a_- +\tilde b_-e^{-z\tau c} =0$. The latter asymptotic formula 
 for $\phi(t)$  is however incompatible with (\ref{40}).  
\hfill \ter
\end{remark}
%
%
Next, the behaviour of a bistable wavefront at $-\infty$ is described in the following proposition: 
\begin{lem}
\label{2af*-}   Let the hypotheses {\bf (B),  (U$^*$)} be satisfied.  If $\phi(t)$ is a bistable wavefront, then there exists a maximal interval $(-\infty, m)$ such that $\phi'(t)>0,\ \phi''(t) >0$ for all $t < m$. Moreover, $\phi(m) \geq e_2$ and,   for some appropriate $t_1\in \R$ and small $\epsilon >0$, 
\begin{equation} \label{as22+}
\phi(t+t_1) = e_1+e^{\mu_1 t} + O(e^{(\mu_1+\epsilon) t}),  \quad \phi'(t+t_1) = \mu_1e^{\mu_1 t} + O(e^{(\mu_1+\epsilon) t}),  \ t \to -\infty. 
\end{equation}
\end{lem}
{\it Proof. } 
Set  $\phi_1(t)= \phi(t)-e_1$ and 
$$
a_j(t) : = \int_0^1g_j(e_1+s\phi_1(t), e_1+s\phi_1(t-c\tau))ds. 
$$
Clearly,  $\phi_1(-\infty)=0, \ \phi_1'(-\infty)=0$ and $a_1(-\infty)= a_+, a_2(-\infty)=b_+$ so that, 
in view of the properties of $\chi_+(z)$ established in Lemma \ref{lc2},  the differential equation for $\phi_1(t)$,
$$
\phi''(t)-c\phi'(t)+a_1(t)\phi(t)+a_2(t)\phi(t-c\tau)=0,
$$
is exponentially dichotomic at $-\infty$.  Moreover, this equation has one-dimensional unstable space
which asymptotically converges to one-dimensional unstable space of the limit equation 
$$
\phi''(t)-c\phi'(t)+a_+\phi(t)+b_+\phi(t-c\tau)=0, 
$$
see \cite[Lemma 4.3]{HL}. This means that $\phi_1'(t) = (\mu_1+o(1))\phi_1(t), \ t \to -\infty,$ and therefore for some $C \not=0$, it holds 
$
\phi_1(t) = C\exp(\mu_1t(1+o(1))), \ t \to -\infty. 
$
If we suppose that $C <0$ than $\phi_1'(t) < 0$ on some maximal interval $(-\infty, s)$ where  $s$ is such that $\phi(s) <  e_1, \ \phi(s)  < \phi(s-c\tau) <0, \ \phi''(s) \geq 0,$ $\phi'(s) =0$.  Consequently, 
$g(\phi(s), \phi(s-c\tau)) \leq 0$,  in contradiction with  $g(\phi(s), \phi(s-c\tau)) > g(\phi(s), \phi(s))>0$. 
Hence, $C>0$ and  $\phi'(t) > 0$ on some maximal interval $(-\infty, r)$.   Suppose that $r$ is finite and 
$\phi(r) < e_2$. Since, in addition,   $\phi(r)  > \phi(r-c\tau) > e_1, \ \phi''(r) \leq 0,$ $\phi'(r) =0$, we obtain that  
$g(\phi(r), \phi(r-c\tau)) \geq 0$,  in contradiction with  $g(\phi(r), \phi(r-c\tau)) <  g(\phi(r), \phi(r))<0$.  The same argument shows that the case $\phi''(r') = 0,$ $\phi'(r') >0$ for some $r' < r$ is not possible as well. Finally, we note that the formulas (\ref{as22+}) is a refinement  of the representation $
\phi_1(t) = C\exp(\mu_1t(1+o(1))), \ t \to -\infty. 
$ Since $
a_1(t) = a_++ O(e^{\nu t}),  \quad a_2(t) =b_+ + O(e^{\nu t}), \ t \to -\infty, 
$ for some positive $\nu$, they can be deduced  from  \cite[Proposition 7.2 ]{MP}, cf. the  proof of Lemma \ref{2af+}. 
 \hfill \ter
 
In the remainder of this section, we assume that $(\tau,c) \in \frak{ D}$ and that the bistable wavefront $\phi$ is monotone.   
After fixing some $\delta \in (0, \mu_1)\cap (0, -\lambda_2)$ such that 
$(1+\gamma)(\mu_2+\delta) < \mu_2$, 
we will consider ${\frak F}_c$ as a linear operator defined
on $C^1_\delta$ and taking its values in $C_\delta$, where 
\[\hspace{-20mm}
C_\delta = \{y= (y_1,y_2) \in C(\R, \R^2): |y|_\delta: = \sup_{s \leq
0}e^{-(\mu_1- \delta)s} |y(s)|+  \sup_{s \geq
0}e^{-(\lambda_{2} + \delta)s} |y(s)| < \infty  \},
\]
\[
C^1_\delta = \{y \in C_\delta: y, y' \in C_\delta, \
|y|_{1,\delta}: =  |y|_\delta +  |y' |_\delta < +\infty
\}. \]
The following result is an immediate consequence of Lemmas \ref{2aff}, \ref{2af*-}. 
  \begin{cor}\label{pmu} Let the hypotheses {\bf (B),  (U$^*$)} be satisfied.  If $\phi(t)$ is a monotone bistable wavefront, then 
$(\phi',\phi'') \in C_\delta$ and 
$
a(t) = a_++ O(e^{\gamma \mu_1 t}),  \quad b(t) = b_++ O(e^{\gamma \mu_1 t}), \ t \to -\infty. 
$
\end{cor}
By repeating the proof of Proposition \ref{P1+} with $(p,q)= (\mu_1-\delta, \lambda_2+\delta)$, we  conclude that 
\begin{pro} \label{P1} ${\frak F}_c:C^1_\delta  \to C_\delta$ is Fredholm operator of index ind ${\frak F}_c=0$. Moreover, ${\frak F}_c$ has one-dimensional  kernel $N({\frak F}_c) =< (\phi'(t), \phi''(t))>$. Thus range $R({\frak F}_c)$ of ${\frak F}_c$ has codimension 1 and therefore 
$$
R({\frak F}_c)= \{F=(f_1,f_2) \in C_\delta: \int_\R \left(f_1(s)v_*(s)+ f_2(s)w_*(s)\right)ds =0\}, 
$$
where $z(t)= (v_*(t),w_*(t))$ is the unique (up to a constant multiple) non-zero solution of  (\ref{asy*}) satisfying inequalities
$
|z(t)| \leq Ke^{-(\lambda_2+\delta)t}, \ t \geq 0; \quad |z(t)| \leq Ke^{-(\mu_1-\delta) t}, \ t \leq 0. 
$
\end{pro}
\begin{remark} \label{RR1}  The asymptotic estimates of $z(t)$ given in Proposition \ref{P1} can be easily improved till 
$$
z(t) =O(e^{-\lambda_1 t}), \ t \geq 0; \quad z(t) = O(e^{-\mu_2 t}), \ t \leq 0. 
$$
For instance, let us prove the first of these formulas. 
Indeed, by Lemma \ref{nt} we can assume that $w_*(s)>0$ for some $s >0$. Set $w(t)= w_*(-t)$ then $w'(t)=v_*(-t)+cw_*(-t)$ and thus 
$$
|(w(t),w'(t))| \leq K_1e^{(\lambda_2+\delta)t}, \ t \leq 0; \quad |(w(t),w'(t))| \leq K_1e^{(\mu_1-\delta) t}, \ t \geq 0. 
$$
Since the characteristic function of the limiting equation for (\ref{rev}) at $-\infty$ is $\chi_-(z)$, \cite[Proposition 7.2 ]{MP} together with Lemma \ref{sms} yield, with for some small $\varepsilon >0$,  the following representation  (possibly, after an appropriate translation of $w(t)$)
$$
w(t) = e^{\lambda_1 t}  + O( e^{(\lambda_1+\varepsilon) t}),\quad w'(t) = \lambda_1 e^{\lambda_1 t}  + O( e^{(\lambda_1+\varepsilon) t}),\quad  t \to -\infty. $$
Therefore 
 $(w_*(t), w_*'(t))  =  (e^{-\lambda_1 t},  -\lambda_1 e^{-\lambda_1 t})  + O( e^{-(\lambda_1+\varepsilon) t}), \ t \to +\infty, $ so that 
 $z(t) =O(e^{-\lambda_1 t}), \ t \geq 0$.  
 
 In addition, we obtain that  $w'_*(t) <0$ for all sufficiently large $t$. Let $d$ be the rightmost critical point of $w_*(t)$. Then $w''_*(d) \leq 0, \ w_*(d) >0, \ w_*(d+c\tau)>0$,  so that equation $w''_*(d)+ cw_*'(d)+ a(d)w_*(d) + b(d+c\tau)w_*(d+c\tau) = 0$ implies that  $d<0$ and $w_*''(t) > 0$ for all $t \geq 0$. 
 \end{remark}

\begin{lem} \label{Lp+}   Let the hypotheses {\bf (B),  (U$^*$)} be satisfied.  Then solution $w_*(t)$ is positive for $t\geq 0$ and non-negative for $t \leq 0:\  w_*(t) \geq 0 , \ t \in \R_-$. 
\end{lem} 
{\it Proof. }  As we have already established  in Remark \ref{RR1},  $w_*(t)>0$ for all $t \geq 0$.  Suppose for a moment that $w_*(t)$ takes negative values on $(-\infty,0)$.  Then there are 
a function $f \in Y_\delta$ and a real number $T>0$ with the following properties
\begin{enumerate}
\item[(i)]  $f(t)=0, \ t \geq 0,$ and $f(t) <0, \ t <0;$ 
\item[(ii)]  $f(t) = {\mathcal D}(-te^{\mu_1 t}) =  -e^{\mu_1 t}\left(\chi_+'(\mu_1)  +o(1))\right) <0, \ t \leq -T,$ (Corollary \ref{pmu} is used here); 
\item[(iii)]  $\int_\R f(t)w_*(t)dt =0$. 
\end{enumerate}
Then, by Corollary \ref{Co2},  the inhomogeneous equation 
\begin{equation}\label{test}
\psi''(t)-c \psi'(t) + a(t)\psi(t) + b(t)\psi(t-c\tau) = f(t)
\end{equation}
has a solution $\psi_* \in X_\delta$.  
Since $f(t) =0$ for $t\geq 0$, and $\psi_*(t)$ is bounded, we conclude that 
$\psi_*(t)/\phi'(t)$ converges to a finite limit as $t\to +\infty$.  
 
 On the other hand,  due to our definition of  $f(t)$ for negative $t$, we obtain that the function $q(t):=\psi_*(t)+te^{\mu_1 t}, \ t \leq - T,$ satisfies the homogeneous equation  ${\mathcal D}q(t)=0, \ t \leq - T$. 
Since $q(t), q'(t)$ have an exponential rate of convergence to $0$ at $-\infty$, we can conclude that $q(t) = B\phi'(t)$  for some finite $B$.  This shows that 
the solution 
$
\psi(t,\xi) = \psi_*(t) + \xi \phi'(t)
$
of (\ref{test}) is positive at $-\infty$ for every real $\xi$.  In this way, 
$
\psi(t,\xi) >0, \  t \in \R, 
$
for all large $\xi >0$.  Set now $$\xi_*: = \inf\{\xi: \psi(t,\xi) >0, \  t \in \R\}. $$ Clearly, $\xi_*$ is finite and 
$\psi(t,\xi) \geq 0, \ t \in \R,$ if and only if $\xi \geq \xi_*$.  Next, since $\psi(t,\xi_*)$ can not have positive maxima on $\R_+$, we obtain  that either (A) $\psi(t,\xi_*)>0, \ \psi'(t,\xi_*)<0$ for $t\geq 0$ 
or (B) $\psi(t,\xi_*)\equiv 0$ on $\R_+$. 

In the case (B), (\ref{test}) implies  that $\psi(s,\xi_*)=0, \ s \in  [-c\tau,0]$, so that  $C^2$-smooth function $\psi(t,\xi_*)$ satisfies  the following algebraic equation 
$$
b(s)\psi(s-c\tau,\xi_*) = f(s), \ s \in  [-c\tau,0].
$$
However, this is not possible because $f(s) <0$ and $b(s) >0, \psi(s-c\tau,\xi_*) \geq 0$ for all $s<0$. 

Now, in the case (A), we have that $\psi(t,\xi_*) >0$ for all $t\in \R$.  Indeed, if 
 $\psi(t_*,\xi_*)=0$ for some $t_* <0$ then 
$\psi''(t_*,\xi_*)\geq \psi'(t_*,\xi_*)=0$  contradicting to (\ref{test}) at $t_*$ (since $f(t_*)<0$, $b(t_*) >0$ and $\psi(t_*-c\tau,\xi_*)\geq 0$). Next, $\psi(t,\xi_*) $
satisfies the differential equation 
\begin{equation}\label{Laa}
\psi''(t)-c \psi'(t) + a_-\psi(t) + b_-\psi(t-c\tau) = n(t), \quad t \in \R, 
\end{equation}
where $n(t)=m(t)+f(t),$
$$
m(t): =  [g_1(e_3,e_3)- g_1(\phi(t),\phi(t-c\tau))]\psi(t,\xi_*) + [g_2(e_3,e_3)- g_2(\phi(t),\phi(t-c\tau))]\psi(t-c\tau,\xi_*)$$
is such that, for some small $\delta_0>0$, it holds 
$$
m(t) \leq 0, \ t \in \R; \  n(t) <0, \ t < 0; \  n(t) = o(1), \ t \to -\infty; \ n(t) = O(e^{(\lambda_2- \delta_0)t}), \ t \to +\infty.
$$
Note that the non-positivity of $m(t)$ follows from the sub-tangency assumption of  {\bf (U$^*$)}. Applying the bilateral 
Laplace transform to (\ref{Laa}) (similarly as it was done in the proof of Lemma \ref{2aff}), we find that, for some  $r\in (0,\delta_0)$, it holds
$$
\psi(t,\xi_*)= \mbox{Res}_{z=\lambda_2} \frac{e^{zt}\tilde n(z)}{\chi_-(z)} + O(e^{(\lambda_2-r) t}), \quad t \to +\infty.
$$
Here $\tilde n(z) = \int_\R e^{-zs}n(s)ds,\  \Re z \in (\lambda_2-\delta_0, 0), $ is the bilateral Laplace transform of $n(t)$. 
A simple calculation shows that 
$$
\mbox{Res}_{z=\lambda_2} \frac{e^{zt}\tilde n(z)}{\chi_-(z)}= \frac{e^{\lambda_2t}\tilde n(\lambda_2)}{\chi'_-(\lambda_2)} = Ae^{\lambda_2t}, \ A:=  \frac{ \int_\R e^{-\lambda_2 s}n(s)ds}{\chi'_-(\lambda_2)}>0, \quad \mbox{if \ } \lambda_3 < \lambda_2; 
$$
$$
\mbox{Res}_{z=\lambda_2} \frac{e^{zt}\tilde n(z)}{\chi_-(z)}=  (Bt +D)e^{\lambda_2t}, \ B:=  \frac{ 2\int_\R e^{-\lambda_2 s}n(s)ds}{\chi''_-(\lambda_2)}>0, \quad \mbox{if \ } \lambda_3 = \lambda_2.  
$$
The described asymptotic behaviour of $\psi(t,\xi_*)>0$ at $\pm \infty$ implies that  
$\psi(t,\xi_*-\epsilon)>0$, $t \in \R$, for all small $\epsilon >0$.  However, this contradicts the definition of $\xi_*$. 
Hence, the non-negativity of $w_*(t)$ on $\R_-$ is proved.  
\hfill \ter
\begin{remark} If, in addition to {\bf (B),  (U$^*$)}, we assume that $g_1(u,v) <0$ for all $(u,v)$ satisfying $u \geq v,$ $u \geq \kappa$, then $w_*(t) >0$ for all $t \in \R$. Indeed,  in such a case,  $a(t)<0$ for all $t \geq -c\tau$.  By arguing as in the last paragraph of Remark \ref{RR1}, this allows to conclude  that $w_*''(t)>0$ for all $t >-c\tau$.  
Let now $m$ be the leftmost point for which the inequality  $w_*(t) >0,\ t \in (m, +\infty)$ holds. Clearly, $m < -c\tau$. 
If $m$ is finite, then $w_*(m)=w_*'(m)= 0 \leq w_*''(m)$. Therefore, in view of equation 
 ${\mathcal D}^*w_*(t)=0$, we find that $b(m+c\tau)w_*(m+c\tau) \leq 0$, contradicting to the fact 
 that $b(m+c\tau) >0, \ w_*(m+c\tau) >0$. 
\end{remark}
\begin{remark} \label{W2} As in  Remark \ref{W1} and Lemma \ref{2affW},    the assumption $(\tau, c) \in \frak{ D}(\tilde a_-,\tilde b_-) \subset \frak{ D}(a_-,b_-) $ can be used instead of the sub-tangency condition
of  Lemma \ref{Lp+}.  Indeed, similarly to  the proof of Lemma \ref{2affW}, it suffices  to replace $(a_-,b_-)$ with $(\tilde a_-,\tilde b_-)$ in formula (\ref{Laa}), and, assuming that $\lambda_3<\lambda_2,$
$$
\psi(t,\xi_*)= A e^{\lambda_3 t } + o(e^{\lambda_3 t} ), \quad \psi'(t,\xi_*) = A \lambda_3e^{\lambda_3 t} + o(e^{\lambda_3 t}), \ t \to +\infty, 
$$
obtain the conflicting representation 
$
\psi(t,\xi_*) = (Pt+Q)e^{\tilde \lambda_2t}(1+o(1)), \ t \to +\infty, \quad |P|+|Q| >0,  
$
with $\tilde \lambda_2\in (\lambda_3, \lambda_2)$ being the biggest negative root of the equation $z^2-cz+\tilde a_- +\tilde b_-e^{-z\tau c} =0$.   \hfill \ter
\end{remark}
\section{Proofs of Theorems \ref{main1a} and \ref{main1b}.}\label{ThT}
In equation (\ref{twe2ng}),   it is convenient to use new independent parameters 
$c, h = c\tau$ instead of $c>0, \tau \geq 0$. Then (\ref{twe2ng}) takes the form 
\begin{eqnarray} \label{twe2nh} \hspace{-15mm}
\phi''(t) - c\phi'(t) + g(\phi(t), \phi(t-h)) =0,  \quad t \in \R, \quad \phi(-\infty)=e_1, \   \phi(+\infty)=e_3.  \end{eqnarray}
\subsection{Local boundedness of the functions $c(h)$ and $c(\tau)$.}
In the coordinates $(c,h)$, the critical curve $c=clin(\tau)$ and the domain  ${\frak D}(a_-,b_-)$ have different shapes described in the following proposition.  Recall  that ${\frak D}(a_-,b_-)$ is defined  as the set of non-negative parameters for which  $\chi_-(z), \ c >0,$ has exactly three real zeros (counting multiplicity). 
\begin{lem} \label{cash}  Set $h_*=\theta(a_-,b_-)>0$. Then there exists a continuous function $c^{\frak E}: \R_+ \to \R_+$, with the properties $c^{\frak E}(h)=0, \ h \in [0,h_*];$  $c^{\frak E}(h)>0, \ h >h_*,$ and $\lim_{h \to \infty} c^{\frak E}(h)/h = 1/\tau_\#$,  such that 
$$ \frak{ D}(a_-,b_-) = \{(h,c): h \geq 0, \  c \geq c^{\frak E} (h)\} \subset \R_+^2. $$
\end{lem}
{\it Proof. }  Since $a_- <0$, it suffices to consider equation (\ref{cL}) for a fixed $h'\geq 0$. Since $\frak{A}(c,h')$ is strictly increasing to $+\infty$ with respect to $c\geq 0$ and  $\frak{B}(c,h')$ is decreasing with respect to $c\geq 0$, this equation have a unique positive  solution $c^{\frak E} (h')$ if and only if $\frak{A}(0,h')\leq \frak{B}(0,h')$. Now, it can be easily verified that the equation $\frak{A}(0,h)=  \frak{B}(0,h)$  has a unique root $h_*= \theta(a_-,b_-)>0$. The computation of the limit $\lim_{h \to \infty} c^{\frak E}(h)/h$ is immediate from equation (\ref{cL}) and the definition of $\tau_\#$ given in Remark \ref{R1*}. 
\hfill \ter 

Next, we show that the velocities of bistable wavefronts are uniformly bounded with respect to $h$ taken from a compact subset of $\R_+$:   
\begin{lem} \label{bau} Suppose that hypothesis  {\bf (B)}  is satisfied. Suppose further that, for each pair $(c_j, h_j)$ of parameters $h_j \in [0,h'],\ c_j>0,\ j \in \N$, problem (\ref{twe2nh}) has  a  monotone  solution $\phi_j: {\mathbb R} \to [e_1,e_3]$.  Then there exists $K=K(h')>0$ such that $c_j \leq K(h'), \ j \in \N$.  
\end{lem}
{\it Proof.} Indeed, suppose that $\epsilon_j = 1/c_j \to 0$. After realising the change of variables $\phi_j(t) = \psi_j(\epsilon_j t)$ and setting $G_j(t)= \psi(t) +g(\psi(t), \psi(t-\epsilon_j h_j)) $, we find that $\psi_j(t)$ satisfies the equation 
\begin{eqnarray} \label{twe2am} \hspace{-15mm}
\epsilon_j^2\psi''(t) - \psi'(t) -\psi(t) = -G_j(t), \ \psi(-\infty)=e_1, \   \psi(+\infty)=e_3.  \end{eqnarray}
Equation (\ref{twe2am}) is translation invariant and therefore we can suppose that $\psi_j(0)= (e_1+e_2)/2$.  
Since $\psi_j(t)$ is a bounded solution of (\ref{twe2am}), it satisfies the integral equation 
\begin{equation}\label{iie} \hspace{0mm}
\psi_j(t) = \frac{1}{\sqrt{1+4\epsilon_j^2}}\left\{\int_{-\infty}^te^{z_j^- (t-s)}G_j(s)ds +
\int_t^{+\infty}e^{z_j^+ (t-s)}G_j(s)ds \right\}.
\end{equation}
where  $z_j^- < 0 < z_j^+$ denote the roots $\epsilon_j^2
z^2 -z -1 =0$.  Clearly, $z_j^- \to -1, \ z_j^+ \to +\infty$.  Differentiating (\ref{iie}), we get 
\begin{equation}\label{iied} \hspace{-15mm}
\psi_j'(t) = \frac{1}{\sqrt{1+4\epsilon_j^2}}\left\{z_j^-\int_{-\infty}^te^{z_j^- (t-s)}G_j(s)ds +z_j^+
\int_t^{+\infty}e^{z_j^+ (t-s)}G_j(s)ds \right\}.
\end{equation}
From (\ref{iied}) we deduce the uniform boundedness of $\psi'_j$: 
$$
|\psi'_j(t)| \leq 2\max\{u+|g(u,v)|, \ u,v \in [e_1,e_3]\}.   
$$
Thus we can find a subsequence $\psi_{j_k}(t)$ of $\psi_j(t)$ which converges, uniformly on compact subsets of $\R$, to  some continuous monotone   function $\psi_*: \R \to [e_1,e_3]$ such that $\psi_*(0)= (e_1+e_2)/2,$ $\psi_*(t) \leq \psi_*(0)$ for 
$t \leq 0$. Invoking the Lebesgue dominated convergence theorem, we find that $\psi_*(t)$ satisfies the integral equation 
$$
\psi_*(t) = \int_{-\infty}^te^{-(t-s)}\left(\psi_*(s) +g(\psi_*(s), \psi_*(s))\right)ds. 
$$
In this way, 
$$
\psi_*'(t)=g(\psi_*(t), \psi_*(t)), \quad  \psi_*(0)= (e_1+e_2)/2, \ \psi_*(t) \leq \psi_*(0) \ \mbox{for all}\  t\leq 0. 
$$
However, due to the bistability of $g(x,x)$ the latter situation is not possible. 
\hfill \ter

A similar result also holds for equation   (\ref{twe2ng}): 
\begin{lem} \label{C5}  Suppose that either the hypothesis  {\bf (U)} or {\bf (U$^*$)}   is satisfied. Then for each $A$ there exists $K(A)>0$ such that $c(\tau) \in (0, K(A)]$ for each  monotone  bistable wavefront   $u=\phi(x+c(\tau)t)$  of equation (\ref{17ng}) considered with $\tau \in [0,A]$. \end{lem}
{\it Proof. }   
We have to prove that  the function $c(\tau)$ is bounded on $[0,A]$. We can argue as in the proof of Lemma  \ref{bau} with the following difference in our reasoning: now we should admit the possibility that the sequence $\tau_j := \epsilon_j h_j \in [0,A]$ can posses a 
subsequence (we will keep the same notation $\tau_j$ for it) converging to a positive limit  $\hat \tau \in [0,A]$. 
Similarly, we will establish the existence of a continuous monotone   function $\psi_*: \R \to [e_1,e_3]$ such that $\psi_*(0)= (e_1+e_2)/2,$ $\psi_*(t) \leq \psi_*(0)$, $t \leq 0$, and  
\begin{equation}\label{gps}
\psi_*'(t)=g(\psi_*(t), \psi_*(t- \hat \tau)). 
\end{equation}
Monotonicity and boundedness of $\psi_*(t)$ also implies that $\psi_*(-\infty) =e_1$ and that $\psi_*(+\infty) \in \{e_2,e_3\}$. 

In particular, there exists $T'$ such that $\psi_*(t) \in (e_1,\kappa]$ for all $t \leq T'$. Therefore, if {\bf (U)} is assumed then  $\psi_*'(t) =g(\psi_*(t), \psi_*(t- \hat \tau))< 0, \ t \leq T'$.  Since  $\psi_*(t)$ is monotone increasing non-constant function, this leads to a contradiction.   
 
On the other hand, if {\bf (U$^*$)} is assumed  then the characteristic equation for linearisation of equation 
(\ref{gps}) 
around the equilibrium $e_1$  is of the form $\lambda +|a_+|= b_+e^{-\lambda \hat \tau}$ with $|a_+|> b_+>0$.  Clearly, all roots of this equation have negative real parts so that the steady state $e_1$ of equation (\ref{gps}) is uniformly asymptotically stable. However, this  is not possible due to  the existence of the solution $\psi_*(t)$  belonging to the unstable manifold of the equilibrium $e_1$. 
 
All the above said proves that  the set $\{c(\tau): \tau \in [0,A]\}$ is bounded. 
\hfill \ter 
\subsection{Local continuation of wavefronts under assumption  {\bf (U)}.}\label{S52}
Assume  {\bf (U)} and  suppose that, given $\tau_0\geq 0$, equation  (\ref{twe2ng})  has a monotone bistable wavefront  $u(t,x)= \phi_0(t+c_0t),$ $c_0>0$.   By Lemma \ref{Lp}, the solution $w_*(t)$ of equation ${\mathcal D}^*w(t) =0$ is non-negative on some maximal interval $[T,+\infty)$ with $T \in [-\infty,0]$.  For our considerations in this section, the case when $T$ is a finite number is much more difficult than the case $T=-\infty$.  Therefore, in what follows, we  assume  that $T \in (-\infty,0]$ (so that $w_*(T)=0$). If $T=-\infty$ then our subsequent arguments simplify with  correctors   $\psi_\varepsilon, S_\varepsilon$ (which are defined below the next lemma) taken identically zero:  
 $\psi_\varepsilon(t) = S_\varepsilon(t)=0$ for all $t \in \R$ (in Section \ref{S53}, these simplifications appear explicitly).  In particular, the following result is needed only when $T \in \R$: 
 
 \begin{lem}  Set $h_0 = c_0\tau_0$ and let $\tilde X_\delta$ and $\tilde Y_\delta$ denote the Banach spaces obtained from $X_\delta$ and $Y_\delta$ by restricting the domain of functions in $X_\delta, Y_\delta$ from $\R$ to $(-\infty,T]$. Then 
for $c,h$ close to $c_0, h_0$,  equation  (\ref{twe2nh}) has a family of solutions $\psi_*(t,c,h), \ t \leq T+h,$ with the following properties:  

\noindent 1) $\psi_*(T,c,h)=\phi_0(T)$ and $\psi_*(t,c,h) >0,$ $ \psi_*'(t,c,h) >0$ for all $t\leq T$; 
\ \  2) $\psi_*(t,h_0,c_0)= \phi_0(t);$ 
 
\noindent 3) $\psi_*(t,c,h)$ depends $C^1$-smoothly on $c,h$ and $(\psi_*)_c=D_c\psi_* :=\partial \psi_*(\cdot,c,h)/\partial c   \in \tilde X_\delta$. 
In particular,  $(\psi'_*)_c = ((\psi_*)_c)'$,  $ ((\psi'_*)_c)'= ((\psi_*)_c)'' \in  \tilde Y_\delta$.  
 \end{lem} 
  {\it Proof. } 
We can consider solution $\psi_*(t,c,h)$ as a perturbation of $\phi_0(t)$: 
$$
\psi_*(t,c,h) = \phi_0(t) +  \zeta (t, c, h), \quad \zeta \in \tilde X_\delta.$$  
Then the equation for $\zeta$ is  
$
{\mathcal D}_0\zeta(t) = N(\zeta, c, h)(t), 
$
where 
\begin{equation}\label{Cd0}
{\mathcal D}_0\zeta(t):=\zeta''(t)-c_0 \zeta'(t) + a(t)\zeta(t) + b(t)\zeta(t-h_0), 
\end{equation}
$$
N(\zeta, c, h) = (c-c_0)(\phi_0'(t)+\zeta'(t)) + a(t)\zeta(t)+b(t)\zeta(t-h_0) +
$$
$$
+g(\phi_0(t),\phi_0(t-h_0)) -  g(\phi_0(t)+\zeta(t),\phi_0(t-h)+\zeta(t-h)), \quad N(0, c_0, h_0) =0.
$$ 
An auxiliary technical result given below, Lemma \ref{dpl}, implies that $N: \tilde X_\delta \times (0,\infty)\times [0,\infty) \to \tilde Y_\delta$ is 
continuously differentiable and $D_\zeta N(0,c_0,h_0)=0$. On the other hand, as it was shown in Claim I 
of Lemma \ref{Lp},  continuous  linear 
operator ${\mathcal D}_0: \tilde X_\delta \to \tilde Y_\delta$ has one-dimensional kernel: dim Ker ${\mathcal D}_0 =1$. 
We claim that, in addition,  ${\mathcal D}_0$  is a surjective operator.   Indeed, take some  $f \in  \tilde Y_\delta$ and  set Set  $f_1(t):=f(t)e^{-(\lambda_1-\delta)t}$. Then consider inhomogeneous equation 
${\mathcal D}_0 u = f$.  The change of variables $u(t) = e^{(\lambda_1-\delta)t}v(t)$ transforms it into
$$
v''(t) - (c_0-2(\lambda_1-\delta))v'(t) + ((\lambda_1-\delta)^2-c_0(\lambda_1-\delta)+a(t))v(t) + e^{-(\lambda_1-\delta)h_0}b(t)v(t-h_0) = f_1(t). 
$$
By our assumptions on $f$ and $\delta$, the function $f_1(t)$ is bounded and 
the limit equation for the latter equation at $-\infty$, 
$$
v''(t) - (c_0-2(\lambda_1-\delta))v'(t) + ((\lambda_1-\delta)^2-c_0(\lambda_1-\delta)+a_-)v(t) + e^{-(\lambda_1-\delta)h_0}b_-v(t-h_0) =0, 
$$
is exponentially dichotomic on $\R$. Then the well known  results from the exponential dichotomy theory (e.g.,  see Lemmas 3.2 and 4.3 in \cite{HL})  show that 
the homogeneous equation 
$$
v''(t) - (c_0-2(\lambda_1-\delta))v'(t) + ((\lambda_1-\delta)^2-c_0(\lambda_1-\delta)+a(t))v(t) + e^{-(\lambda_1-\delta)h_0}b(t)v(t-h_0) = 0 
$$
possesses an exponential dichotomy on $(-\infty,T]$  so that the above considered 
 inhomogeneous equation  has at least one bounded solution $v_*(t), \ t \leq T$, with $v_*'(t), v_*''(t)$ which are also bounded on $\R_-$. 
It is clear that $u_*(t) = e^{(\lambda_1-\delta)t}v_*(t) \in \tilde X_\delta$ and ${\mathcal D}_0 u_* = f$. 

The smoothness properties  of operator $N$ and the Fredholm property of ${\mathcal D}_0$ allow to realize a standard Lyapunov-Schmidt reduction in the equation ${\mathcal D_0}\zeta(t) = N(\zeta, c, h)(t)$.  The details of this procedure (used in more complex situation) are described in Lemma \ref{L12} below.  This method allows to establish 
the existence of one-parametric family of functions $\zeta = \zeta(t, h, c, a)$ depending $C^1$-smoothly on parameters $(h,c,a)$ close to $(h_0,c_0,0)$ and such  that 
$\zeta(t,h_0,c_0, 0) =0$ for all $t \leq T$ and 
$
 \phi_0(t) +  \zeta (t, c, h, a)
$
solves equation  (\ref{twe2nh}) for each fixed $(c,h,a)$. Note that the dimension $1$ of parameter $a$ corresponds to  dim Ker ${\mathcal D}_0 =1$. 
We can  reintroduce this parameter in a more usual way by fixing  
$a=0$ and considering  the family  of shifted solutions $\phi_0(t+s)+\zeta (t+s, c, h, 0)$ of equation (\ref{twe2nh}). 

Finally, consider the equation 
$\phi_0(T+s)+\zeta (T+s, c, h, 0) = \phi_0(T)$. Clearly $(s,h,c)= (0,h_0,c_0)$ is a solution of this equation while $\phi_0'(T)>0$. Therefore, 
in view of the implicit function theorem, there exists $C^1-$smooth solution $s = s(c,h)$ of this equation satisfying equality $s(h_0,c_0)=0$.   
We obtain the required family $\psi_*$ by setting 
$$
\psi_*(t,c,h) = \phi_0(t+s(c,h)) +  \zeta (t+s(c,h), c, h,0). 
$$
Observe also that the monotonicity properties of  $\psi_*(t,c,h)$ are assured by Lemma \ref{2aaf} and 
$$
\frac{\partial \psi_*(\cdot,c,h)}{\partial c} = \left( \phi_0'(\cdot+s(c,h)) +  \zeta' (\cdot+s(c,h), c, h,0)\right)s_1(c,h) +  \zeta_2(\cdot+s(c,h), c, h,0) \in  \tilde X_\delta.   \ter
$$ 

Next, using $C^\infty-$smooth non-increasing function $S_\varepsilon(t)$ such that $S_\varepsilon(t) =1$ for $t\leq T$ and $S_\varepsilon(t) =0$ for $t\geq T+\varepsilon$, we will define  the `corrector' 
$\psi_\varepsilon(t,c,h) = (\psi_*(t,c,h)-\phi_0(t))S_\varepsilon(t)$. Clearly, $\partial \psi_*(T,c,h)/\partial c =0$ and 
\begin{eqnarray*}\label{psn}
\psi_\varepsilon(t,c,h)= \left\{\begin{array}{cc} 
\psi_*(t,c,h)-\phi_0(t),& t \leq T, \\    
0, & t \geq  T+ \varepsilon.\end{array}\right. 
\end{eqnarray*}
Also  $\psi_\varepsilon(t,c_0,h_0) \equiv 0$. 
We will look for a monotone solution $\phi(t,c,h), \ t \in \R,$ of (\ref{twe2nh}) in the form 
$$
\phi(t,c,h) = \phi_0(t) + \psi_\varepsilon(t,c,h)+ \zeta (t, c, h), 
$$
where $\zeta \in X_\delta$.  
Then the equation for $\zeta$ is  
$
{\mathcal D}_0\zeta(t) = N_\varepsilon(\zeta, c, h)(t), 
$
where 
$
{\mathcal D}_0$
is given by (\ref{Cd0}) and  $
N_\varepsilon(\zeta, c, h)(t) = $ $$
(c-c_0)(\phi_0'(t)+\zeta'(t)) +c((\psi_*(t) - \phi_0(t))S_\varepsilon(t))' - ((\psi_*(t) - \phi_0(t))S_\varepsilon(t))'' + a(t)\zeta(t)+
$$
$$
b(t)\zeta(t-h_0) +
g(\phi_0(t),\phi_0(t-h_0)) -  g(\phi_0(t)+\psi_\varepsilon(t)+\zeta(t),\phi_0(t-h)+\psi_\varepsilon(t-h)+\zeta(t-h)), \ t \in \R.  
$$
Since $\psi_*(t,c,h) = \phi_0(t)+ \psi_\varepsilon(t,c,h),\ t \leq T,$ solves equation (\ref{twe2nh})  for all $t \leq T$, it is easy to find that 
$$
N_\varepsilon(\zeta, c, h)(t) = (c-c_0)\zeta'(t) + a(t)\zeta(t)+b(t)\zeta(t-h_0) +
$$
$$
+g(\psi_*(t,c,h),\psi_*(t-h,c,h))-  g(\psi_*(t,c,h)+\zeta(t),\psi_*(t-h,c,h)+\zeta(t-h)), \quad t\leq T. 
$$
$N_\varepsilon$ has the following smoothness properties: 
 \begin{lem}  \label{L19} There exist neighborhoods $\mathcal{O}(c_0)$ and $\mathcal{O}(h_0)$ of the points $c_0,h_0$ such that 
 function $N_\varepsilon: X_\delta \times \mathcal{O}(c_0) \times \mathcal{O}(h_0) \to Y_\delta$, $
N_\varepsilon(0, c_0,h_0)=0,$  is 
continuously differentiable, with continuous  partial derivatives given by
$$
D_cN_\varepsilon(\zeta, c, h)(t) = \phi_0'(t) +\zeta'(t)-(\phi_0(t)S_\varepsilon(t))' + S_\varepsilon(t)D_cg(\psi_*(t),\psi_*(t-h))  - (\psi_*)_c(t)S''_\varepsilon(t)+
$$
$$
(\psi_*(t)+c(\psi_*)_c(t)-2(\psi_*'(t))_c) S'_\varepsilon(t)  -  D_cg(\phi_0(t)+\psi_\varepsilon(t)+\zeta(t),\phi_0(t-h)+\psi_\varepsilon(t-h)+\zeta(t-h)), \ t \in \R;
$$
$$
D_cN_\varepsilon(\zeta, c, h)(t) = \zeta'(t) 
+D_cg(\psi_*(t,c,h),\psi_*(t-h,c,h))-  $$
$$D_cg(\psi_*(t,c,h)+\zeta(t),\psi_*(t-h,c,h)+\zeta(t-h)), \ t\leq T. 
$$
$$D_hN_\varepsilon(\zeta, c, h) = S_\varepsilon(t) D_hg(\psi_*(t,c,h), \psi_*(t-h,c,h)) + S_\varepsilon'(t) (c(\psi_*)_h(t) - 2(\psi_*')_h(t)) - $$
$$
-(\psi_{*})_h(t)S''_\varepsilon(t) - D_hg(\phi(t,c,h), \phi(t-h,c,h));
$$
$$D_\zeta N_\varepsilon(\zeta, c, h)w(t) = (c-c_0)w'(t) + a(t)w(t)+b(t)w(t-h_0)- $$
$$
g_1(\phi_0(t)+\psi_\varepsilon(t)+\zeta(t),\phi_0(t-h)+\psi_\varepsilon(t-h)+\zeta(t-h))w(t) - 
$$
$$
g_2(\phi_0(t)+\psi_\varepsilon(t)+\zeta(t),\phi_0(t-h)+\psi_\varepsilon(t-h)+\zeta(t-h))w(t-h). 
$$
In particular, $D_\zeta N_\varepsilon(0, c_0, h_0)=0,$
\begin{eqnarray*}\label{DNC}
D_c N_\varepsilon(0, c_0, h_0)(t) = \left\{\begin{array}{cc} 
0,& t \leq T, \\    
\phi_0'(t) + R_\varepsilon(t), & t \geq  T,\\
\phi_0'(t),& t \geq T+\varepsilon,   
\end{array}\right. 
\end{eqnarray*}
where 
$
R_\varepsilon(t)  = S_\varepsilon(t)D_cg(\psi_*(t),\psi_*(t-h_0))  + 
(c_0(\psi_*)_c(t)-2(\psi_*'(t))_c) S'_\varepsilon(t) - 
(\psi_*)_c(t)S''_\varepsilon(t) +$  

\noindent $-D_cg(\phi_0(t)+\psi_\varepsilon(t),\phi_0(t-h_0)+\psi_\varepsilon(t-h_0)) -\phi'_0(t)S_\varepsilon(t) . 
$ 
\end{lem}
{\it Proof. }  The proof of this lemma is based on routine straightforward calculations. Some of them (concerning  $D_\zeta N_\varepsilon$) are given below, in the proof of a similar technical assertion, Lemma \ref{dpl}. Here it is convenient to use the relation  
$\mathcal{T}:=c((\psi_*(t) - \phi_0(t))S_\varepsilon(t))' - ((\psi_*(t) - \phi_0(t))S_\varepsilon(t))'' =
-c (\phi_0(t)S_\varepsilon(t))' + $ $$(\phi_0(t)S_\varepsilon(t))'' +S_\varepsilon(t)g(\psi_*(t,c,h),\psi_*(t-h,c,h)) +(c\psi_*(t)-2\psi'_*(t))S'_\varepsilon(t) -\psi_*(t)S''_\varepsilon(t).
$$
Recall that $\psi_*(t,\cdot,h): \mathcal{O}(c_0) \to \tilde X_\delta$ depends $C^1$-continuously on $c\in \mathcal{O}(c_0)$ while all terms of $\mathcal{T}$ belong to the space $Y_\delta$ since $S_\varepsilon(t) =0$ for all $t \geq T+\varepsilon$.
\hfill \ter
\begin{cor} \label{NZ} It holds that 
$$
\lim_{\varepsilon \to 0^+}\int_{-\infty}^{+\infty}w_*(t)D_c N_\varepsilon(0, c_0, h_0)(t)dt =  \lim_{\varepsilon \to 0^+} \int_T^{T+\varepsilon}w_*(t)R_\varepsilon(t)dt+  \int_T^{+\infty}w_*(t)\phi'_0(t)dt= $$
$= \int_T^{+\infty}w_*(t)\phi'_0(t)dt >0, 
$
so that $\int_{-\infty}^{+\infty}w_*(t)D_c N_\varepsilon(0, c_0, h_0)(t)dt >0$ for all small positive $\varepsilon$.  
\end{cor}
{\it Proof. }  Indeed, by integrating by parts and using the boundary conditions $w_*(T)=(\psi_*)_c(T)=0$, we find that $\int_T^{T+\varepsilon}w_*(t)R_\varepsilon(t)dt  =$
$$
\int_T^{T+\varepsilon}w_*(t)\left\{S_\varepsilon(t)D_cg(\psi_*(t),\psi_*(t-h_0)) -D_cg(\phi_0(t)+\psi_\varepsilon(t),\phi_0(t-h_0)+\psi_\varepsilon(t-h_0))\right\}dt$$
$$-\int_T^{T+\varepsilon}S_\varepsilon(t)\left( \phi'_0(t)w_*(t)+
 \left[w_*(t)
(c_0(\psi_*)_c(t)-2(\psi_*'(t))_c)\right]' +
(w_*(t)(\psi_*)_c(t))'' \right)
dt =O(\epsilon).  \ter
$$
\begin{lem} \label{L12} Suppose that $\phi_0'(t) >0, c_0 >0,$ and that hypothesis   {\bf (U)}  is satisfied. Then there exist an open neighbourhood ${\mathcal O}$ of $h_0= \tau_0c_0,$ and  $C^1$-smooth  function $c: {\mathcal O} \to (0,+\infty),\ c(h_0) =c_0,$ such that equation  (\ref{twe2nh})  has a continuous 
family $\phi(\cdot, c(h),h)\in \phi_0 + X_\delta,$ $\phi(t, c(h_0),h_0)= \phi_0(t)$, of strictly increasing bistable wavefronts. \end{lem}
{\it Proof. }  Consider the direct sums of Banach spaces $X_\delta= Ker\, {\mathcal D}_0 \oplus W$, $Y_\delta = R ({\mathcal D}_0) \oplus V$, where $Ker\, {\mathcal D}_0 $ is one-dimensional null space of the linear operator ${\mathcal D}_0$ and the range $R({\mathcal D}_0)$ has codimension one,
$$
R({\mathcal D}_0) = \{f\in Y_\delta: \int_\R f(s)w_*(s)ds=0\}, $$
$$V = <y_*> \mbox{for some} \ y_* \in Y_\delta \ \mbox{with}  \ \int_\R w_*(u)y_*(u)du \not =0. 
$$
Let $P: Y_\delta \to Y_\delta$ be the projection on the subspace $R({\mathcal D}_0)$ along $V$,   
$$
Pf(s)= f(s) - \frac{ \int_\R f(u)w_*(u)du}{ \int_\R w_*(u)y_*(u)du}y_*(s). 
$$
Then the equation $
{\mathcal D}_0\zeta = N_\varepsilon(\zeta, c, h), \ \zeta \in X_\delta, 
$
is equivalent to the system 
$$
{\mathcal D}_0\xi = PN_\varepsilon(\xi+u, c, h), \ \xi \in W \subset X_\delta, \ u = k\phi'_0,  \quad  (I-P)N_\varepsilon(\xi+u, c, h)=0
$$
in the sense that $\zeta = \xi +u$ satisfies the former equation if and only if it satisfies the latter system. 
Considering the restriction ${\mathcal D}'= {\mathcal D}_0|_W: W\to R({\mathcal D})$, we know that the operator ${\mathcal D}'$ is invertible and thus the equation ${\mathcal D}'\xi(t) = PN_\varepsilon(\xi+u, c, h)$ can be written as 
$\xi = ({\mathcal D}')^{-1}PN_\varepsilon(\xi+u, c, h)=Q(\xi+u,c,h)$. Since $D_\xi Q(0,c_0,h_0) =0$, the implicit function theorem \cite{AP} shows  that this equation has a $C^1$-continuous family of solutions $\xi= \xi(u,c,h)$ defined in some vicinity of the 
point $(0,c_0,h_0)$, where $\xi(0,c_0,h_0)=0$. We still need to prove that for appropriate parameters $(c,h)$ close to $(c_0,h_0)$ the equation 
$$
(I-P)N_\varepsilon(\xi(u,c,h)+u, c, h)=0
$$
is satisfied. Simplifying, we can take $u=0$. Since  for all small $\varepsilon>0$, in view of Corollary \ref{NZ},  $$D_c(I-P)N_\varepsilon(\xi(0,c,h_0), c, h_0)|_{c=c_0}= (I-P)D_cN_\varepsilon(0,c_0,h_0)=  $$ $$\frac{ \int_\R D_c N_\varepsilon (0, c_0, h_0)(u) w_*(u)du}{ \int_\R w_*(u)y_*(u)du}y_*(s)\not= 0,$$
we conclude that there exists a $C^1-$continuous solution $c=c(h), \ c(h_0)=c_0$, $h \in {\mathcal O}$, of the equation $(I-P)N_\varepsilon(\xi(0,c,h), c, h)=0$.  To finalise the proof of the lemma, we have to establish the monotonicity  of  bistable waves 
$$
\phi(t,c(h),h) := \phi_0(t) + \psi_\varepsilon(t,c,h) +\xi (0, c(h), h)(t), \ h \in {\mathcal O}. 
$$
First, note that Lemmas \ref{cot}, \ref{2aaf} imply that 
$\phi(\cdot,c(h),h): {\mathbb R}\to (e_1,e_3)$ for all $ h\in {\mathcal O}$. Moreover,  each  $\phi(t,c(h),h)$ is strictly monotone in $t$ on some maximal interval $(-\infty,s_h)$,  where $\phi(s_h,c(h),h) > \kappa$, and is also strictly monotone at $+\infty$, see Lemma \ref{2af+}.  In fact, we prove below that the following asymptotic formula holds at $+\infty$: 
\begin{equation} \label{apo}\hspace{-25mm}
\phi'(t,c(h),h) = q_1(h)e^{\mu_2 t} + r(t,h), \quad \mbox{where} \ |r(t,h)| \leq Ke^{(\mu_2-\delta')t}, \ t \geq 0,  \ h \in  {\mathcal O}' \subset {\mathcal O}, 
\end{equation}
where $K\geq 1, \delta' >0$ does not depend on $h$, and $q_1(h)$ is a positive continuous function defined on some smaller neighbourhood ${\mathcal O}'$ of $h_0$.  It follows from (\ref{apo}) that $\phi'(t,c(h),h) >0$ for 
all $h \in  {\mathcal O}'$ and $t > t_* := (\delta')^{-1} \sup\{\ln (K/q_1(h)),  h \in  {\mathcal O}'\}$.  Since $\phi'(t,c(h),h) $ converges to $\phi'_0(t)>0$ (as $h \to h_0$) uniformly on compact subsets of $\mathbb R$, we may conclude  that $\phi'(t,c(h),h)  >0$ for all $t \in \mathbb R$ once $h$ is sufficiently close to $h_0$. 

To prove (\ref{apo}), we apply the bilateral Laplace transform to the differential equation  for $\psi(t) : =e_3 - \phi(t,c(h),h)$: 
\begin{equation}\label{pps}
\psi''(t)-c\psi'(t) + a_+\psi(t) + b_+\psi(t-h)=  d(t,h), 
\end{equation}
where $d(t,h) = -(g(e_3,e_3)- g (e_3-\psi(t), e_3-\psi(t-h)) - a_+\psi(t) - b_+\psi(t-h))= O(e^{(\mu_2(h_0)+\delta)(1+\gamma)t}),$ $ t \to +\infty$. Importantly, since the function 
$\xi(0,c(h),h):  {\mathcal O} \to X_\delta$ is continuous, the latter $O$ relation for $d(t,h)$ is satisfied  uniformly with respect to   $h$ from compact subsets of  ${\mathcal O}$. 
Thus,  for some small positive $r>0$ satisfying $\mu_2(h)-r > (\mu_2(h_0)+\delta)(1+\gamma)$
and $t >0$, we have that 
$$
\psi(t)= \mbox{Res}_{z=\mu_2} \frac{e^{zt}\tilde d(z,h)}{\chi_+(z)} + \frac{1}{2\pi i} \int_{\Re z = \mu_2-r} \frac{e^{tz}\tilde d(z,h)}{\chi_+(z)}dz = \alpha(h) e^{\mu_2t} + B(t,h),$$ with $ B(t,h) := \beta(t,h) e^{(\mu_2-r)t}.
$
Here $\tilde d(z,h) = \int_\R e^{-zs}d(s,h)ds,\  \Re z \in ((1+\gamma)(\mu_2(h_0)+\delta), 0), $ is the bilateral Laplace transform of $d(t,h)$. In view of the Lebesgue dominated convergence theorem and  the uniform (with respect to $h$ from compact subsets of ${\mathcal O}$) exponential estimate  $d(t,h) = O(e^{(\mu_2(h_0)+\delta)(1+\gamma)t}),$ $ t \to +\infty$, the transform $\tilde d(z,h)$ depends continuously on $z,h$ and is uniformly bounded on the vertical line $\{\Re z = \mu_2-r\}$. 

Consequently, in view of Lemma \ref{2af+}, there exists some small neighbourhood ${\mathcal O}' \subset {\mathcal O}$ of $h_0$ such that continuous functions $\alpha(t), \ \beta(t,h)$ satisfy the estimates
$$
\alpha(h) = \frac{\tilde d(\mu_2(h),h)}{\chi'_+(\mu_2(h))} >0,  \quad |\beta(t,h)|\leq  \frac{1}{	2\pi} \int_{\mathbb R} \frac{|\tilde d(\mu_2(h)-r +is,h)|}{|\chi_+(\mu_2(h)-r +is)|}ds \leq B_0, \ h \in {\mathcal O}',
$$
where  $B_0= B_0({\mathcal O}')$ is some  positive constant.  Finally, integrating $(\ref{pps})$ on $(t,+\infty)$, we find that
$$
\psi'(t) = c\psi(t) +\int_t^{+\infty}(a_+\psi(s) + b_+\psi(s-h)-  d(s,h))ds =\mu_2(h)\alpha(h)e^{\mu_2(h)t}+ R(t,h), 
$$ 
where 
$$
R(t,h) = cB(t,h) +\int_t^{+\infty}(a_+B(s,h) + b_+B(s-h,h)-  d(s,h))ds
$$
satisfies $|R(t,h)| \leq De^{(\mu_2-r)t}, t > 0$, $h \in  {\mathcal O}'$, with $D$ not depending  on $h$. 
\hfill \ter
\subsection{Local continuation of wavefronts under hypothesis   {\bf (U$^*$)}.} \label{S53}
When {\bf (U$^*$)} is assumed instead of {\bf (U)}, the local continuation of wavefronts is somewhat  easier to prove. The main reason of this  is  the non-negativity of   solution $w_*(t)$ of the adjoint equation.   Indeed,  at the beginning of  Subsection  \ref{S52}, we mentioned that  the proofs in this subsection simplify  when  $w_*(t) \geq 0,$ $ t \in \R$ (i.e. when $T=-\infty$).     Therefore in  Subsection  \ref{S52} we narrowed our attention to more complex case of finite  $T$.  In the present subsection we show how the Lyapunov-Schmidt reduction works for $T=-\infty$. 
 
Hence,  suppose that, given $\tau_0\geq 0$, equation  (\ref{17ng})  has  a monotone bistable wavefront  $u(t,x)= \phi_0(t+c_0t),$ $c_0>0$.  For $c,h$ close to $c_0, h_0= c_0\tau_0$, we will look for a monotone solution $\phi(t,c,h)$ of (\ref{twe2nh}) in the form 
$$
\phi(t,c,h) = \phi_0(t) + \zeta (t, c, h), 
$$
where $\zeta \in X_\delta$.  
Then the equation for $\zeta$ is  
$
{\mathcal D}_0\zeta(t) = N(\zeta, c, h), 
$
where 
$
{\mathcal D}_0\zeta, \ 
N(\zeta, c, h)$ are defined in (\ref{Cd0}). 
Next result can be regarded as somewhat simplified version of Lemma \ref{L19}:
 \begin{lem}  \label{dpl} Function $N: X_\delta \times (0,+\infty)\times [0,+\infty) \to Y_\delta$ is 
continuously differentiable, with continuous  partial derivatives given by
$D_cN(\zeta,c,h) = \phi_0'(t)+ \zeta'(t)$, 
$$D_hN(\zeta, c, h) = g_2(\phi_0(t)+\zeta(t),\phi_0(t-h)+\zeta(t-h))(\phi_0'(t-h)+\zeta'(t-h));$$
$$D_\zeta N(\zeta, c, h)w(t) = (c-c_0)w'(t) + a(t)w(t)+b(t)w(t-h_0)
$$
$$
-  g_1(\phi_0(t)+\zeta(t),\phi_0(t-h)+\zeta(t-h))w(t) - g_2(\phi_0(t)+\zeta(t),\phi_0(t-h)+\zeta(t-h))w(t-h). 
$$
In particular, 
$
N(0, c_0,h_0)=0, \quad D_\zeta N(0, c_0, h_0)=0, \quad D_c N(0, c_0, h_0)=\phi_0'(t),$ $$ D_hN(0, c_0, h_0) = g_2(\phi_0(t),\phi_0(t-h_0))\phi_0'(t-h_0).
$$
\end{lem}
{\it Proof. }  Clearly, it suffices to check the validity of the conclusions of Lemma \ref{dpl} only for the nonlinear part $N_1$ of $N$. Here 
 $
N_1(\zeta, h) = g(\phi_0(t),\phi_0(t-h_0)) -  g(\phi_0(t)+\zeta(t),\phi_0(t-h)+\zeta(t-h)), 
$
and below we will give details of computations only for  more difficult derivative $D_\zeta N_1$, the other derivatives being similar. To abbreviate, we use the notation $f_h(t)=f(t-h)$.  First, we find that 
$
\Delta:=$ $$N_1(\zeta+w, h) - N_1(\zeta, h) + g_1(\phi_0(t)+\zeta(t),\phi_h(t)+\zeta_h(t))w(t) + g_2(\phi_0(t)+\zeta(t),\phi_h(t)+\zeta_h(t))w_h(t)= 
$$
$$
\int_0^1\left(g_1(\phi_0(t)+\zeta(t),\phi_h(t)+\zeta_h(t))- g_1(\phi_0(t)+\zeta(t)+sw(t),\phi_h(t)+\zeta_h(t)+sw_h(t)\right)ds \, w(t)+
$$
$$
\int_0^1\left(g_2(\phi_0(t)+\zeta(t),\phi_h(t)+\zeta_h(t))- g_2(\phi_0(t)+\zeta(t)+sw(t),\phi_h(t)+\zeta_h(t)+sw_h(t)\right)ds \, w_h(t). 
$$
Therefore, for every $r>0$ there exists  $K_r$ such that for all  $w$ such that $|w|_{\infty} \leq r$, it holds 
$$
|\Delta| \leq K_r\left(|w(t)|^{1+\gamma}+ |w(t-h)|^\gamma |w(t)|+ |w(t)|^\gamma |w(t-h)|+ |w(t-h)|^{1+\gamma}\right).    
$$
The latter implies that $|\Delta|_{Y_\delta} \leq K_r'\left(|w|^{1+\gamma}_{Y_\delta}\right)$  for some $K_r' \geq K_r$ and  all $w$ such that $|w|_{Y_\delta} \leq r$.  This proves that the Fr\'echet derivative $D_\zeta N_1$ exists and is given by 
$D_\zeta N_1(\zeta,  h)w(t) =$ $$ -  g_1(\phi_0(t)+\zeta(t),\phi_0(t-h)+\zeta(t-h))w(t) - g_2(\phi_0(t)+\zeta(t),\phi_0(t-h)+\zeta(t-h))w(t-h). 
$$
Next,  it can proved similarly that $D_\zeta N_1(\zeta,  h)$ is locally H\"older continuous function in view of the estimate
$$
\|D_\zeta N_1(\zeta_1,  h_1)- D_\zeta N_1(\zeta,  h)\| \leq K\left(|\zeta_1-\zeta|_\infty^{\gamma}+ [|\zeta'|_\infty+|\phi'_0|_\infty]^{\gamma}|h-h_1|^{\gamma}+ |h-h_1|\right). \quad \ter 
$$
\begin{lem} \label{L12p} Suppose that $\phi_0'(t) >0, c_0 >0$ and that  hypothesis  {\bf (U$^*$)} with $c_0 < clin(h_0/c_0)$ is satisfied. Then there exist an open neighbourhood ${\mathcal O}$ of $h_0= \tau_0c_0,$ and  $C^1$-smooth  function $c: {\mathcal O} \to (0,+\infty),\ c(h_0) =c_0,$ such that equation  (\ref{twe2nh})  has a continuous 
family $\phi(\cdot, c(h),h)\in \phi_0 + X_\delta,\  h \in {\mathcal O},$ $\phi(t, c(h_0),h_0)= \phi_0(t)$, of strictly increasing bistable wavefronts. If {\bf (U$^*$)} holds with $c_0 =clin (h_0/c_0)$, the same conclusion, possibly except for the  strict monotonicity property of $\phi(\cdot, c(h),h)$ at $+\infty$,  is true. 
\end{lem}
{\it Proof. }  Taking the non-negative solution $w_*(t)$ defined in Section \ref{S4}, we consider the  Banach spaces $W$, $V$ and the projector 
 $P: Y_\delta \to Y_\delta$ defined in the first paragraphs of the proof of Lemma \ref{L12}. 
Then the equation $
{\mathcal D}_0\zeta(t) = N(\zeta, c, h), \ \zeta \in X_\delta, 
$
is equivalent to the system 
$$
{\mathcal D}_0\xi(t) = PN(\xi+u, c, h), \ \xi \in W \subset X_\delta, \ u = k\phi'_0,  \quad  (I-P)N(\xi+u, c, h)=0.
$$
Considering the restriction ${\mathcal D}'= {\mathcal D}_0|_W: W\to R({\mathcal D})$, we know that the operator ${\mathcal D}'$ is invertible and thus the equation ${\mathcal D}'\xi(t) = PN(\xi+u, c, h)$ can be written as 
$\xi = ({\mathcal D}')^{-1}PN(\xi+u, c, h)=Q(\xi+u,c,h)$. Since $D_\xi Q(0,c_0,h_0) =0$,  this equation has a $C^1$-continuous family of solutions $\xi= \xi(u,c,h)$ defined in some vicinity of the 
point $(0,c_0,h_0)$, where $\xi(0,c_0,h_0)=0$. We have to prove that for appropriate parameters $(c,h)$ close to $(c_0,h_0)$ the equation 
$$
(I-P)N(\xi(u,c,h)+u, c, h)=0
$$
is satisfied. It suffices to take $u=0$. Since $$D_c(I-P)N(\xi(0,c,h_0), c, h_0)|_{c=c_0}= (I-P)D_cN(0,c_0,h_0)=  \frac{ \int_\R \phi'_0(u)w_*(u)du}{ \int_\R w_*(u)y_*(u)du}y_*(s)\not= 0,$$
we conclude that there exists a $C^1-$continuous solution $c=c(h), \ c(h_0)=c_0$, $h \in {\mathcal O}$, of the equation $(I-P)N(\xi(0,c,h), c, h)=0$.  

For $c_0 > c^{\frak E} (h_0)$, we will prove now  the monotonicity  of the obtained  bistable waves 
$$
\phi(t,c(h),h) := \phi_0(t) + \xi (0, c(h), h)(t), \ h \in {\mathcal O}.
$$
  The restriction $c_0 > c^{\frak E} (h_0)$ implies  that $\chi_-(z)$ has exactly three different real zeros, $\lambda_3 < \lambda_2 <0< \lambda_1$. 
By Lemma \ref{2af*-}, $\phi(t,c(h),h)$ is strictly monotone in $t$ on some maximal interval $(-\infty,r_h)$,  where $\phi(r_h,c(h),h) \geq  e_2$. Therefore, to complete the proof of Lemma \ref{L12}, it suffices to prove  the following asymptotic formula (which is similar to (\ref{apo})):
\begin{equation} \label{apor}\hspace{-25mm}
\phi'(t,c(h),h) = q_2(h)e^{\lambda_2 t} + r_2(t,h), \ \mbox{where} \ |r_2(t,h)| \leq K_2e^{(\lambda_2-\delta'')t}, \ t \geq 0,  \ h \in  {\mathcal O}'' \subset {\mathcal O}, 
\end{equation}
where $K_2\geq 1, \delta'' >0$ does not depend on $h$, and $q_2(h)$ is a positive continuous function defined on some small neighbourhood ${\mathcal O}''$ of $h_0$. Indeed, once (\ref{apor}) is established, we can argue as in the paragraph below formula  (\ref{apo}).  

Now, in order to prove (\ref{apor}), we will apply the bilateral Laplace transform to the differential equation  for $\psi(t) : =e_3 - \phi(t,c(h),h)$: 
\begin{equation}\label{pps2}
\psi''(t)-c\psi'(t) + a_-\psi(t) + b_-\psi(t-h)=  d_*(t,h), 
\end{equation}
where $d_*(t,h) = -(g(e_3,e_3)- g (e_3-\psi(t), e_3-\psi(t-h)) - a_-\psi(t) - b_-\psi(t-h)), \ t \in \R$. Clearly,  in view of the sub-tangency restriction imposed in   {\bf (U$^*$)}, 
$$d_*(t,h_0) = -(g(e_3,e_3)- g (e_3-\phi_0(t), e_3-\phi_0(t-h)) - a_-\phi_0(t) - b_-\phi_0(t-h_0)) \leq 0, \ t \in \R.$$
Also, $d_*(t,h_0) \not\equiv 0$ on $\R$ and  
$d_*(t,h) = O(e^{(\lambda_2(h_0)+\delta)(1+\gamma)t}),$ $ t \to +\infty$, with $O$ relation being satisfied  uniformly with respect to   $h \in  {\mathcal O}''$. 
Thus,  for some small positive $r'>0$ satisfying $\lambda_2(h)-r' > (\lambda_2(h_0)+\delta)(1+\gamma)$
and $t >0$, we have that $\tilde d_*(\lambda_2(h_0),h_0) >0$, 
$$
\psi(t)= \mbox{Res}_{z=\lambda_2} \frac{e^{zt}\tilde d_*(z,h)}{\chi_-(z)} + \frac{1}{2\pi i} \int_{\Re z = \lambda_2-r'} \frac{e^{tz}\tilde d_*(z,h)}{\chi_-(z)}dz = \alpha(h) e^{\lambda_2t} + B(t,h), $$ with  $B(t,h) := \beta(t,h) e^{(\lambda_2-r')t}.
$
Here $\tilde d_*(z,h) = \int_\R e^{-zs}d_*(s,h)ds,\  \Re z \in ((\lambda_2(h_0)+\delta)(1+\gamma), 0), $ is the bilateral Laplace transform of $d_*(t,h)$. In view of the Lebesgue dominated convergence theorem and  the uniform exponential estimate  $d_*(t,h) = O(e^{(\lambda_2(h_0)+\delta)(1+\gamma)t}),$ $ t \to +\infty$,  $h \in  {\mathcal O}''$, $\tilde d_*(z,h)$ depends continuously on $z,h$. 

As consequence, there exists some small neighbourhood ${\mathcal O}''' \subset {\mathcal O}$ of $h_0$ such that continuous functions $\alpha(h), \ \beta(t,h)$ satisfy the estimates
$$
\alpha(h) = \frac{\tilde d_*(\lambda_2(h),h)}{\chi'_-(\lambda_2(h))} >0,  \quad |\beta(t,h)|\leq  \frac{1}{	2\pi} \int_{\mathbb R} \frac{|\tilde d_*(\lambda_2(h)-r' +is,h)|}{|\chi_-(\lambda_2(h)-r' +is)|}ds \leq B_0, \ h \in {\mathcal O}''',
$$
where  $B_0= B_0({\mathcal O}''')$ is some  positive constant.  Finally, integrating $(\ref{pps})$ on $(t,+\infty)$, we find that
$$
\psi'(t) = c\psi(t) +\int_t^{+\infty}(a_-\psi(s) + b_-\psi(s-h)-  d_*(s,h))ds =\lambda_2(h)\alpha(h)e^{\lambda_2(h)t}+ R(t,h), 
$$ 
where 
$$
R(t,h) = cB(t,h) +\int_t^{+\infty}(a_-B(s,h) + b_-B(s-h,h)-  d_*(s,h))ds
$$
satisfies $|R(t,h)| \leq De^{(\mu_2-r)t}, t > 0$, $h \in  {\mathcal O}'''$, with  $D$ not depending on $h$. 
\hfill \ter
\begin{remark} \label{W3} Note that  the representation (\ref{apor}) remains valid if 
 the sub-tangency condition
of  Lemma \ref{Lp+} is replaced with the assumption $(\tau, c) \in \frak{ D}(\tilde a_-,\tilde b_-) \subset \frak{ D}(a_-,b_-) $.  Indeed, in such a case, Lemma \ref{2affW} assures  that $q_2(h_0) \not =0$.  
\end{remark}
\subsection{Global continuation of wavefronts.}
Now we can complete the proof of Theorems \ref{main1a}, \ref{main1b}.  Consider the family $\frak{F}$ of all continuous functions 
$c_\alpha:[0,h_\alpha) \to (0,+\infty), \ \alpha \in A,$ such that for every $h=c\tau \in [0,h_\alpha)$ equation (\ref{twe2ng}) has a bistable monotone wavefront propagating with the velocity $c_\alpha(h)$ and $c_\alpha(0)=c_0$ where $c_0$ is the speed of the unique bistable monotone front of the non-delayed equation. Lemmas \ref{L12}, \ref{L12p} show that $\frak{F}$  is a non-empty set, $A \not= \emptyset$.  We will introduce a partial order $\prec$ in $\frak{F}$ in the following way: $(c_\alpha, h_\alpha) \prec (c_\beta, h_\beta)$ if $h_\beta \geq h_\alpha$ and $c_\alpha(h)= c_\beta(h)$  for all $h\in[0,h_\alpha)$.  Clearly, we can apply the Zorn lemma to the family $(\frak{F}, \ \prec)$, let 
$c^*:[0,h^*) \to (0,+\infty)$ be the maximal element. Note that if {\bf (U$^*$)} is assumed then  the graph $\frak{G}$ of the curve $c^*$ belongs to the domain 
$\frak{ D}(a_-,b_-)$.  

Suppose first that $h^* = +\infty$, then $\sup_{h \geq 0} h/c^*(h)=+\infty$ since otherwise there exists a bounded sequence of delays 
$\tau_j = h_j/c^*(h_j), \ h_j \to +\infty, $ such that $c^*(h_j) = h_j/\tau_j \to + \infty$, contradicting to the conclusion of Lemma \ref{C5}.  
In view of the intermediate value theorem,    this implies that for each 
$$\tau  \in [0, +\infty) = [\min_{h \geq 0} h/c^*(h), \sup_{h \geq 0} h/c^*(h))$$
there exists at least one monotone bistable wavefront.  

Now, if $h^* < +\infty$, then,  by Lemma \ref{bau}, 
$c^*$ is a bounded function on $[0,h^*)$. Let the interval $[p,q]$ (it can happen that $p=q$) denote  the set of all partial limits of $c^*(h)$ as $h \to h^*-$.  Set $r = (p+q)/2$ and suppose that $r >0 $ (and, in addition,  $r > c^{\frak E} (h^*)$  if {\bf (U$^*$)} is assumed). Then  there exists a sequence $h_j \to h^*$ and $c_j=c^*(h_j)$ such that $c_j \to r$.  The sequence of profiles $\phi_j$ of wavefronts $\phi_j(x+c_jt)$ is uniformly bounded and equicontinuous on $\R$.  We can also assume that $\phi_j(0)= (e_1+e_2)/2$ for every $j$. Therefore we can find a
subsequence of $\phi_j$ (we will use the same notation $\phi_j$ for it) converging to some non-decreasing function 
$\phi_*$ such that $\phi_*(0) = (e_1+e_2)/2$ and $\phi_*(t) \in [e_1,e_3], \ t \in \R$.  It is easy to see that $\phi_*$ satisfies the differential equation
$$
\phi''(t) - r\phi'(t) + g(\phi(t), \phi(t-h^*)) =0, 
$$ 
and $\phi_*(\pm\infty ) \in \{e_1,e_2,e_3\}$.  Actually, since $\phi_*$ is non-decreasing and $\phi_*(0) <e_2$, we obtain that $\phi_*(-\infty)=e_1$.  Considering the possibility $\phi_*(+\infty)= e_2$, we find that $\phi_*(t) \in [e_1,e_2]$ for all $t\in \R$ and therefore $\phi_*''(t) \geq 0, \ t \in \R$. Clearly, this contradicts to the convergence of 
$\phi_*(t)$ at $+\infty$. Thus $\phi_*(t)$ is a strictly monotone bistable wavefront of the above limiting delay  differential equation.  Consequently, we can apply either  Lemma \ref{L12} or Lemma  \ref{L12p} for parameters $c=r, h= h^*$ and conclude that 
there exists a smooth function $c=c(h)$ with $h$ from some open neighbourhood  $\mathcal U$ of $h^*$ and a family of bistable wavefronts $\phi_{c(h)},\  \phi_r=\phi_*$ for  all $h \in \mathcal U$.
 In addition, if {\bf (U)} is assumed, these  wavefronts are monotone. 
 Since $c(h)$ is smooth, we can use this function to extend continuously $c^*$ on the open interval $[0,h^{**})$ strictly bigger than  $[0,h^{*})$. 
This shows  that either $r=0$ or $r= c^{\frak E} (h^*) >0$. 

As we have observed, under conditions of Theorem \ref{main1a}, the only case $r=0$ can happen. In such a case, the graph $\frak{G}$ of $c^*$ connects continuously points $(0,c_0)$ and $(h^*,0)$ so that  for every fixed nonnegative $\tau$ the  line $h= c\tau$ will intersect $\frak{G}$ at least once at some point $(c(\tau)\tau, c(\tau))$. 
This means that if $r=0$ (in particular, this always occurs under conditions of Theorem \ref{main1a}) then  for each fixed $\tau\geq 0$ the original equation has at least one bistable wavefront propagating with the velocity $c(\tau)>0$.  

On the other hand, under conditions of Theorem \ref{main1b},  the situation when 
$c^*(h^*-)=c^{\frak E} (h^*) >0$ can also occur. Then the graph $\frak{G}$ of $c^*$ connects continuously points $(0,c_0)$ and $(h^*,c^* (h^*))$ so that  for every fixed  $\tau \in [0, \tau_*], \ \tau_*  = h^*/c^*(h^*),$ the  line $h= c\tau$ intersects $\frak{G}$ at least once at some point $(c(\tau)\tau, c(\tau))$. 
Thus for each fixed $\tau \in [0, \tau_*]$ the original equation has at least one bistable wavefront propagating with the velocity $c(\tau)>0$.  
Next,  due to the maximality property of $h^*$ there exist arbitrarily small positive $h-h^*$ such that  $(h,c^*(h))\not \in \frak{D}(a_-,b_-)$. Thus for delay $\tau' = h/c^*(h)$ 
(which can be chosen arbitrarily close to $\tau_*$)  there exists a  non-monotone wavefront $\phi$ propagating with the velocity $c^*(h)$ close to $c^{\frak E} (h^*)$. 
Moreover,  since $(h,c^*(h))\not \in \frak{D}(a_-,b_-)$, the leading asymptotic term of $e_3-\phi(t)$ at $+\infty$ is oscillatory,  e.g. see \cite[Lemma 4.6]{GTLMS}. 
Thus $\phi(t)$ is oscillating around $e_3$ at $+\infty$.
\begin{remark} \label{W4} If  the sub-tangency condition of  {\bf (U$^*$)} is not assumed, the above proof remains  true if we consider domain $\frak{ D}(\tilde a_-,\tilde b_-)$ instead of $\frak{ D}(a_-,b_-)$. See Remarks \ref{W1}, \ref{W2}, \ref{W3}.  However, if $\frak{ D}(\tilde a_-,\tilde b_-)\not= \frak{ D}(a_-,b_-)$, then our approach does not allow to extend the curve $\frak{G}$ till the boundary of $\frak{ D}(a_-,b_-)$ making a conclusion about the existence of the oscillating wavefronts. It is worth noting that delayed reaction-diffusion equations can possess non-monotone non-oscillating wavefronts, cf. \cite{IGT}. 
\end{remark}
\section*{Acknowledgements.} \noindent The first author was supported by FONDECYT (Chile) under project 1150480. The second author was partially supported by the Ministry of Education and Science of the Russian Federation (the agreement number 02.a03.21.0008).

\end{document}